\documentclass[psamsfonts]{conm-p-l}
\usepackage{latexsym}
\usepackage{diagrams}
\diagramstyle[scriptlabels,small]
\usepackage{amssymb}

\newarrow{Equal}{}{=}{}{=}{}

\def\lto{{\longrightarrow}}
\def\into{{\hookrightarrow}}
\def\xto{\xrightarrow}
\def\onto{\twoheadrightarrow}

\newcommand{\fp}{{\mathfrak p}}

\newcommand{\calm}{{\mathcal M}}

\newcommand{\calz}{{\mathcal Z}}

\newcommand{\uEnd}{\underline {\End}}

\newcommand{\oA}{{\overline {A}}}

\newcommand{\oC}{{\overline {C}}}

\newcommand{\oH}{{\overline {H}}}

\newcommand{\oSG}{{\overline {SG}}}

\newcommand{\bbbb}{{\mathbb B}}

\newcommand{\bbbe}{{\mathbb E}}

\newcommand{\bbbz}{{\mathbb Z}}

\newcommand{\bdot}{\bullet}

\newcommand{\vp}{\varphi}

\DeclareMathOperator{\Ann}{Ann}

\DeclareMathOperator{\Coker}{Coker}

\DeclareMathOperator{\Der}{Der}

\DeclareMathOperator{\depth}{depth}

\DeclareMathOperator{\ev}{ev}
\DeclareMathOperator{\End}{ End}

\DeclareMathOperator{\Ext}{Ext}

\DeclareMathOperator{\gldim}{gldim}
\DeclareMathOperator{\grade}{grade}

\DeclareMathOperator{\lgrade}{lgrade}

\DeclareMathOperator{\Hoch}{HH}

\DeclareMathOperator{\Hom}{Hom}
\DeclareMathOperator{\id}{id}

\DeclareMathOperator{\Imm}{Im}

\DeclareMathOperator{\Ker}{Ker}

\DeclareMathOperator{\op}{op}

\DeclareMathOperator{\Spec}{Spec}

\DeclareMathOperator{\Tor}{Tor}
\DeclareMathOperator{\tr}{tr}

\DeclareMathOperator{\Mod}{\ensuremath{ \mathbf{Mod}}}

\theoremstyle{definition}
\newtheorem{defn}{Definition}[section]
\newtheorem{definition}[defn]{Definition}

\newtheorem{rem}[defn]{Remark}
\newtheorem{remark}[defn]{Remark}

\newtheorem{remarks}[defn]{Remarks}

\newtheorem{sit}[defn]{}

\newtheorem{example}[defn]{Example}

\theoremstyle{plain}
\newtheorem{prop}[defn]{Proposition}
\newtheorem{proposition}[defn]{Proposition}
\newtheorem{theorem}[defn]{Theorem}
\newtheorem{lem}[defn]{Lemma}
\newtheorem{lemma}[defn]{Lemma}
\newtheorem{cor}[defn]{Corollary}
\newtheorem{corollary}[defn]{Corollary}
\newtheorem*{Conjecture}{Conjecture}

\theoremstyle{remark}

\begin{document}
\title[Hochschild Cohomology and Idempotents] {Morita
Contexts, Idempotents, and Hochschild Cohomology\\ --- with
Applications to Invariant Rings ---}

\dedicatory{Dedicated to the memory of Peter Slodowy}
\author[R.-O.~Buchweitz]{Ragnar-Olaf~Buchweitz}
\address{Department of Computer and Mathematical Sciences,
University of Toronto at Scarborough, Toronto, Ontario,
Canada M1C 1A4} \email{ragnar@math.utoronto.ca}

 \thanks{The author
gratefully acknowledges partial support from NSERC}

\date{\today} 

\begin{abstract}
We investigate how to compare Hochschild cohomology of
algebras related by a Morita context.  Interpreting a Morita
context as a ring with distinguished idempotent, the key
ingredient for such a comparison is shown to be the grade of
the Morita defect, the quotient of the ring modulo the ideal
generated by the idempotent.  Along the way, we show that
the grade of the stable endomorphism ring as a module over
the endomorphism ring controls vanishing of higher groups of
selfextensions, and explain the relation to various forms of
the Generalized Nakayama Conjecture for Noetherian algebras. 
As applications of our approach we explore to what extent
Hochschild cohomology of an invariant ring coincides with
the invariants of the Hochschild cohomology.
\end{abstract}

\maketitle

{\footnotesize\tableofcontents}

\section*{Introduction}
One of the basic features of Hochschild (co-)homology is its
invariance under derived equivalence, see \cite{Kel}, in
particular, it is invariant under Morita equivalence, a
result originally established by Dennis-Igusa \cite{DIg},
see also \cite{Lod}.  This raises the question how
Hochschild cohomology compares in general for algebras that
are just ingredients of a Morita context.

For example, when a finite group $G$ acts on an algebra $S$,
then the invariant ring $R=S^{G}$ and the skew group algebra
$S\#G$ are related by a Morita context whose other
ingredients are two copies of $S$, one considered as
$(R,S\#G)$--bimodule, the second as $(S\#G, R)$--bimodule. 
It is known, see \cite{Lor, Ste}, how to relate the
Hochschild (co-)\-homology of $S$ to that of $S\#G$, but the
relation of either to that of $R$ seems less well
understood.

However, the analoguous question for tangent (or
Andr\'e-Quillen) cohomology, that is, the cohomology of the
cotangent complex, has been addressed as then one of the
first applications of the theory.  Indeed, M.~Schlessinger
showed in \cite{Schl} that the quotient of a rigid complex
analytic singularity by a finite group is again rigid,
provided the depth along the branch locus of the group
action is at least equal to 3.  His ideas show more
generally that, for $S$ the commutative ring of the
singularity and $R=S^{G}$, the tangent cohomology satisfies
$T^{i}_{R}\cong (T^{i}_{S})^{G}$ for each $i\le g-2$, where
$g$ is the depth along the branch locus.  Our main result in
this direction, \ref{thm:key}, shows this to remain true if
tangent cohomology is replaced by Hochschild cohomology.

Our starting point, in Section 1, is to view a Morita
context as a ring $C$ with a distinguished pair of
complementary idempotents $e, e'=1-e$, relating the rings
$A=eCe$ and $B=e'Ce'$ through the bimodules $M=e'Ce$ and
$N=eCe'$.  We define the associated {\em Morita defects\/}
as $\oC = C/CeC$ and $\oC' = C/Ce'C$, noting that $C$
represents a Morita equivalence between $A$ and $B$ if, and
only if, both defects vanish.

The key notion of the grade of a module and its basic
properties are reviewed in Section 2.  The crucial
observation then is that $C$ is the classical Morita context
defined by a ring $A$ and a module $M$, or what we call an
{\em Auslander context\/}, if, and only if, the grade of the
Morita defect $\oC$ is at least 2.

The natural question as to the significance of that grade in
general has a surprisingly satisfying answer: that grade,
call it $g$, indicates precisely the vanishing of
$\Ext^{i}_{A}(M\oplus A,M\oplus A)$ for $1\le i\le g-2$, as
we show in \ref{thm:grade}.  This makes evident the relation
to various formulations, mainly by M.~Auslander and
I.~Reiten, of the Generalized Nakayama Conjecture that ask,
in essence, for which classes of rings $A$ and classes of
modules $M$ over it, vanishing of those $\Ext$--groups for
each $i\ge 1$ forces $M$ to be projective.  This angle is
explored in Section 3, where we view these conjectures as
attempts to extend Morita theory.  Informally, the key
question becomes, whether, or when, infinite grade of a
Morita defect forces its vanishing.

In Section 4, we finally turn our attention to Hochschild
cohomology.  After introducing the map that compares
Hochschild cohomology of a Morita context with that of one
of its components, we show it to be a homomorphism of graded
algebras.  We added details on this result after the referee
of an earlier version of the paper brought the work of
E.~Green and \O.~Solberg \cite{GSo} to our attention that
deals with this comparison map for Morita contexts under
some additional hypotheses.  The referee also pointed to
various recent papers, such as \cite{Cib,CMRS, GSo, MPl},
that set up exact sequences for Hochschild cohomology of
triangular algebras, a degenerate case of Morita contexts,
but generalizing Happel's original sequence for
one-point-extension algebras \cite{Hap}. The occurring
comparison maps of Hochschild cohomology there are special cases
of the one studied here, as pointed out already in \cite{GSo}.

Our main application of the grade of a Morita defect to
Hochschild cohomology in Section 5 shows coincidence of
those groups for algebras in a Morita context almost up to
the grade.

To apply the results of Section 5 efficiently, the question
that remains is how to determine the grade of a Morita
defect.  We do not know a general answer.  However, for
invariants under a finite group acting on a commutative
noetherian domain, we relate in Section 7 that grade to the
depth of the Noether different of the ring over its subring
of invariants, after setting the stage in Section 6 by
reviewing in general the homological algebra of finite
group actions and explaining the notion of the {\em
Noether\/} or {\em homological\/} different as originally
defined by Auslander-Goldman \cite{AGo} in the
noncommutative setting.

It is not until Section 7 that we deal with commutative
algebra proper, whereas the sections leading up to it should
be seen as advocacy that sometimes excursions into
noncommutative algebra help to shed light on problems in
commutative algebra.

As the reader might easily guess from the references used,
the ideas here are greatly influenced by, if not taken
straight from the work of M.~Auslander and his
collaborators.  Indeed, temptation was great to head this
article ``Homological theory of idempotents'', not only
paying hommage to \cite{APT}, but also retaining the title
of the original talk given at the Joint AMS/SMF meeting in
Lyon, July 2001, of which this is a vastly extended and
detailed version.

I also wish to thank the referee for a careful reading.

\section{Morita Contexts}

We begin with the fairly obvious observation that a ({\em
generalized\/}) {\em Morita context\/} or {\em
pre-equivalence\/} is nothing but a ring with a
distinguished idempotent. This point of view is at least
implicit whenever the data of a Morita context are arranged
as $(2\times 2)$-matrices, see, for example, \cite{MRo}.

\begin{sit}
    Given an idempotent $e=e^{2}$ in a ring $C$, let
    $e'=1-e$ denote its complementary idempotent. The
    resulting Pierce decomposition $C = e'Ce'\oplus
    e'Ce\oplus eCe'\oplus eCe $ into a direct sum of abelian
    groups allows to display $C$ as an algebra of $(2\times
    2)$-matrices,
    $$
    C \cong
    \begin{pmatrix}
        e'Ce'& e'Ce\\
        eCe'&eCe
    \end{pmatrix}\,.
    $$
\end{sit}

\begin{sit}
    The diagonal entries, $A=eCe$ and $B=e'Ce'$, are rings
    with multiplicative identity $1_{A}=e$ and $1_{B}=e'$
    respectively, but they are usually not (unital) subrings
    of $C$. In terms of the rings $A$ and $B$, the component
    $M=e'Ce\subseteq C$ is a left $B$-, right $A$-module,
    whereas the component $N=eCe'\subseteq C$ is a right
    $B$-, left $A$-module. Emphasizing $A, B$ and the
    bimodules $M,N$, the displayed decomposition of $C$
    becomes
    $$
    C \cong
        \begin{pmatrix}
        e'Ce'& e'Ce\\
        eCe'&eCe
        \end{pmatrix} =
        \begin{pmatrix}
         B&M\\
        N&A \end{pmatrix}\,.
    $$
    The multiplication on $C$ induces
    \begin{itemize}
    \item a homomorphism of $B$-bimodules $$
    f:M\otimes_{A}N=e'Ce\otimes_{\displaystyle eCe}eCe'\to
    e'Ce'= B\,,$$
     
    \item a homomorphism of $A$-bimodules $$
    g:N\otimes_{B}M=eCe'\otimes_{\displaystyle e'Ce'}e'Ce\to
    eCe= A\,,$$
    \end{itemize}
    and associativity of the multiplication on $C$ combines
    the associativity conditions required of the Morita
    context $(A,B,M,N,f,g)$. Conversely, arranging a Morita
    context into a $(2\times 2)$-matrix as above defines a
    ring $C$ with the distinguished pair of complementary
    idempotents $e=1_{A}, e' = 1_{B}$. We thus don't
    reinvent the wheel, when we make the
\end{sit}

\begin{definition}
    A {\em Morita context\/}, between rings $A$ and $B$,
    consists of a ring $C$ together with an idempotent 
    $e\in C$, such that $A\cong eCe, B\cong e'Ce'$.
    A morphism of Morita contexts $(C_{1},e_{1})\to 
    (C_{2},e_{2})$
    is a ring homomorphism $\vp:C_{1}\to C_{2}$ such that
    $\vp(e_{1})=e_{2}$.
\end{definition}

Morita contexts were originally introduced to characterize
when two rings have equivalent module categories. The
crucial role of the given idempotents leads us to the
following notion.
\begin{definition}
    The homomorphic images
    $$
    \oC = C/(e) = C/CeC \quad\text{and}\quad \oC'= C/(e') =
    C/Ce'C 
    $$
    are the {\em Morita defects\/} of $C$. When $\oC=0$, we
    call $e$ a {\em Morita idempotent}.
\end{definition}

Note the following simple fact.

\begin{lemma}
    Let $(C,e)$ be a Morita context. A right module $X$ over
    $C$ is annihilated by $CeC$, thus naturally a right
    $\oC$-module, if and only if $Xe=0$. Furthermore, the 
    following two conditions are equivalent:
    \begin{enumerate}
	\item $e$ is a Morita idempotent.
    
	\item A right $C$--module $X$ satisfies $Xe=0$ if
	and only if $X=0$.	
    \end{enumerate}
\end{lemma}

\begin{proof}
    For the first claim, observe that evidently $Xe=0
    \Longleftrightarrow XCeC=0$. The claimed equivalence
    follows from $\oC e = 0$.
\end{proof}

The classical theorem on Morita equivalence is essentially a
consequence of the following, fundamentally important, exact
sequence of $C$-bimodules,
\begin{equation}
\label{eq:1}
\tag{$\Diamond$}
    0\to \Omega_{C/A}\to Ce\otimes_{\displaystyle eCe} 
    eC\xto{\mu_{e}} C\to \oC\to 0\,,  
\end{equation}
where $\mu_{e}$ is induced by the multiplication on $C$, and
$\Omega_{C/A}$ is defined as the kernel of $\mu_{e}$.

\begin{lemma}
\label{lem:fundamental}
    Applying $e'C\otimes_{C} - \otimes_{C}Ce'$ to the
    exact sequence {\em (\ref{eq:1})\/} above results in the
    exact sequence of $B=e'Ce'$-bimodules
\begin{diagram}
    0&\rTo &\Omega_{C/A}&\rTo &e'Ce\otimes_{\displaystyle
    eCe} eCe'&\rTo^{\scriptstyle
    e'C\otimes_{C}\mu_{e}\otimes_{C}Ce'}&
    e'Ce'&\rTo&\oC&\rTo &0\hphantom{\,.}\\
    &&\dEqual&&\dEqual&&\dEqual&&\dEqual\\ 0&\rTo
    &\Omega_{C/A}&\rTo &M\otimes_{A} N&\rTo^{f}& B&\rTo&
    \oC&\rTo &0\,.
\end{diagram}    

In particular, $\oC \cong B/((e)\cap B)$, whence the Morita
defect $\oC$ is a homomorphic image of $B$, isomorphic to
the cokernel of $f$.

The $C$-bimodule $\Omega_{C/A}$ is annihilated by $e$ on
both sides, thus naturally a $\oC$-bimodule.
\end{lemma}

\begin{proof}
    Applying $C\otimes_{C}-\otimes_{C}C$ does not change the
    exact sequence (\ref{eq:1}). Now, as a right $C$-module,
    $C=eC\oplus e'C$, the tensor product $eC\otimes_{C}-$ is
    an exact functor, and $eC\otimes_{C}\mu_{e}$ is clearly
    an isomorphism. In particular, $e\Omega_{C/A}=0$. On the
    other side, $C=Ce\oplus Ce'$ as left $C$--modules, and
    the result follows by combining these observations.
\end{proof}

\begin{rem}
One may show directly, cf.~\cite[Lemma 1.5, Prop.4.6]{APT},
that
$$
\Omega_{C/A} \cong \Tor_{2}^{C}(\oC,\oC)
$$
as $C$--bimodules, and then the last claim in 
\ref{lem:fundamental} becomes completely transparent.    
\end{rem}

\begin{sit}
    Henceforth, without further specification, a ``module''
    over some ring $A$ will mean a unital right module, and
    the category formed by those will be denoted $\Mod A$. 
    Accordingly, we write $\Hom_{A}$ to denote homomorphisms
    of right modules, $\Hom_{A^{\op}}$ to denote those of
    left $A$--modules, where $A^{\op}$ is the opposite ring. 
    Modules over commutative rings are considered symmetric
    bimodules, as usual.
    
    Recall that a module $X$ over a ring $A$ is a {\em
    generator\/} if there exists an $A$-linear epimorphism
    $X^{n}\onto A$, equivalently, $A$ is a direct summand of
    a finite direct sum of copies of $X$. If $X$ is any
    module, then clearly $X\oplus A$ is a generator.
\end{sit}

The following classical result, see \cite[II.Thm.3.4]{Bas},
characterizes a Morita idempotent, and lists its crucial
properties for the associated Morita context. It is 
essentially an immediate consequence of \ref{lem:fundamental}.

\begin{proposition}
\label{prop:Morita}
    The following conditions are equivalent for an idempotent 
    $e\in C$.
\begin{enumerate}
    \item  $e\in C$ is a Morita idempotent, 

    \item  the restricted multiplication map
    $$
    \mu_{e}:Ce\otimes_{\displaystyle eCe}eC\to C
    $$
    is an epimorphism of $C$-bimodules, and in that case, it 
    is an isomorphism.

    \item  the restricted multiplication map 
    $$
    f: N\otimes_{A}M =
    e'Ce\otimes_{\displaystyle eCe}eCe'\to e'Ce'=B
    $$ 
    is an epimorphism of $B$-bimodules, and in that case, it
    is an isomorphism.
\end{enumerate}   
Moreover, if $e\in C$ is a Morita idempotent, then
\begin{enumerate}
    \item[(i)] the $B = e'Ce'$-modules\footnote{violating
    our general convention, here one of them is a right, the
    other a left module, as is obvious from the
    context\ldots} $M=eCe'$ and $N=e'Ce$, as well as $Ce'$
    and $e'C$, are {\em generators\/}.

    \item[(ii)] 
    \label{prop:Moritaii}
    the $A=eCe$-modules $M=eCe'$ and $N=e'Ce$, as well
    as $eC$ and $Ce$, are {\em finite projective\/}.    

    \item[(iii)] the multiplication on $C$, equivalently,
    the adjoints to $g:M\otimes_{B}N\to A$, yield {\em
    isomorphisms\/}
    \begin{align*}
	M=eCe'&\xto{\cong}\Hom_{(eCe)}(e'Ce,eCe)
	=\Hom_{A}(N,A)\\
	&\quad\text{of $(B,A)$-bimodules, and}\\
	N= e'Ce&\xto{\cong}\Hom_{eCe^{\op}}(eCe',eCe) 
	=\Hom_{A^{\op}}(M,A)\\
	&\quad\text{of $(A,B)$-bimodules.}
    \end{align*}

     \item[(iv)] the multiplication on $C$ yields {\em ring
     isomorphisms\/}
\begin{xxalignat}{3}
& &B=e'Ce'\xto{\cong}&\End_{eCe}(e'Ce)=\End_{A}(M)&&\\
    &\hphantom{\square}
	    &B=e'Ce'\xto{\cong} &\End_{(eCe)^{\op}}(eCe')^{\op} =
		    \End_{A^{\op}}(N)^{\op}\,.
		    &&\square
			    \end{xxalignat}
\end{enumerate}
\end{proposition}

Classical Morita theory asserts that $(A,B,M,N,f,g)$
represents a Morita equivalence if and only if $f$ and $g$
are surjective. The latter condition is equivalent to the
vanishing of both Morita defects by the preceding results.
Indeed, we have the following.

\begin{corollary}
\label{cor:ME}
Let $(C,e)$ be a Morita context.
\begin{enumerate}
    \item
   \label{ME:1}
    The rings $A = eCe$ and $C$ are Morita equivalent through
    the bimodules $Ce, eC$, if and only if $\oC=0$.

    \item
    \label{ME:2} The rings $B= e'Ce'$ and $C$ are Morita
    equivalent through the bimodules $Ce', e'C$, if and
    only if $\oC'=0$.
    
    \item 
    \label{ME:3}
    The rings $A$ and $B$ are Morita equivalent via $(C,e)$
    if and only if $\oC =0$ and $\oC'=0$.

\end{enumerate}
\end{corollary}

\begin{proof}
It clearly suffices to prove (\ref{ME:1}). The data
$(A,C,Ce,eC,\mu_{e}, e\mu_{C}e:eC\otimes_{C}Ce\to A)$ form a
Morita context, in which the last map is surjective by
definition of $A$. The map corresponding to $f$ is
$\mu_{e}$, whence the result.
\end{proof}

We cannot resist rephrasing this as well in categorical
terms, see, for example, \cite[5.3, 7.1]{Aus1}, \cite{Gab}, 
as this formulation exhibits clearly the role of the Morita 
defect when comparing module categories.

\begin{theorem}
\label{prop:catME}
Let $(C,e)$ be a Morita context and set $A=eCe$ as before. 
Restriction of scalars embeds $\Mod \oC$ as a
Serre\footnote{This property is equivalent to $\oC$ being
defined by an idempotent ideal, see \cite{Aus1}!}
subcategory into $\Mod C$, and the quotient category $\Mod
C/\Mod \oC$ is equivalent to $\Mod A$ under the functor
$-\otimes_{C}Ce:\Mod C\to \Mod A$.  There is thus an exact
sequence of abelian categories
\begin{diagram}
    0&\rTo &\Mod \oC& \rTo & \Mod C 
    &\rTo^{-\otimes_{C}Ce}&\Mod A&\rTo & 0\,.
\end{diagram}
Moreover, both the inclusion and the projection functor
admit both a left and a right adjoint, and the adjoints,
$-\otimes_{A}eC$ and $\Hom_{A}(Ce,-),$ to $-\otimes_{C}Ce$
are fully faithful.

Clearly, $-\otimes_{C}Ce$ is an equivalence if and only 
if $\Mod \oC = 0$, that is, $\oC = 0$.\qed
\end{theorem}

\section{The Grade of the Stable Endomorphism Ring}

The main theme here are the homological properties of the
Morita defects.  The key invariants that we exploit are
their {\em grades\/} as right or left $C$--modules
respectively.

\begin{definition}
    Recall that for a module $X$ over a ring $A$, its {\em
    grade\/} is given by
    $$
    \grade_{A}X = \inf\{i\ge 0\mid \Ext_{A}^{i}(X,A)\ne
    0\}\quad\in\quad\bbbz\cup\{\infty\}\,.
    $$
    If $X$ happens to be an $A$-bimodule, we will write
    $\lgrade_{A}(X)$ for its grade as left module.
\end{definition}

%

We will use repeatedly the following simple fact.

\begin{lemma}
\label{lem:grade}
    Let $\oA$ be a homomorphic image of the ring $A$. For
    any $\oA$-module $X$ and any finite projective
    $A$--module $P$, one has
    $$
    \Ext^{i}_{A}(X,P) = 0\quad\text{for}\quad 
    i<\grade_{A}\oA\,.
    $$
\end{lemma}

\begin{proof}
    The change-of-rings spectral sequence
    $$
    E^{ij}_{2}= \Ext^{i}_{\oA}(X, \Ext^{j}_{A}(\oA,P))
    \Rightarrow \Ext^{i+j}_{A}(X,P)
    $$
    establishes the claim as $ \Ext^{j}_{A}(\oA,P) = 0$, and
    thus $E^{ij}_{2} = 0$, for $j<\grade_{A}\oA$.
\end{proof}

Coming back to our study of a ring $C$ with distinguished
idempotent $e$, we saw above that $\oC$ is naturally a
homomorphic image of $B=e'Ce'$, whence any $\oC$--module can
be viewed naturally as a $B$--module.

\begin{lemma}
\label{lem:gradechange}
Assume $\oC'=0$, so that $B$ and $C$ are Morita equivalent.
In that case,
$$
\grade_{B}\oC = \grade_{C}\oC\,.
$$
\end{lemma}

\begin{proof}
By \ref{prop:catME}, the functor $-\otimes_{C}Ce'$ is an 
equivalence from $\Mod C$ to $\Mod B$. One has 
$\oC\otimes_{C}Ce' = \oC e' = \oC$, and 
$C\otimes_{C}Ce'\cong e'Ce'\oplus eCe' = B\oplus N$ as 
$B$--modules, with $N$ finite projective by 
\ref{prop:Morita}(ii). Thus,
$$
\Ext^{i}_{C}(\oC,C) \cong \Ext^{i}_{B}(\oC,B\oplus N)\,
$$
and these groups vanish if, and only if,
$\Ext^{i}_{B}(\oC,B)$ vanishes.
\end{proof}

The following example plays a key role.

\begin{example}
    Let $A$ be a ring and $M$ a right $A$-module. The
    projection from $M\oplus A$ onto $A$ defines the
    idempotent $e=(M\oplus A \onto A\into M\oplus A)$ in the
    endomorphism ring $C= C(A,M)= \End_{A}(M\oplus A)$ to
    yield a Morita context between $A$ and $B=\End_{A}(M)$,
    $$
    C(A,M) \cong
    \begin{pmatrix}
    B = \End_{A}(M)& M=\Hom_{A}(A,M)\\
    M^{*} = \Hom_{A}(M,A)& A=\End_{A}(A)
    \end{pmatrix}
    $$
We will call $(C(A,M),e)$ the {\em Auslander context\/}
defined by the pair $(A,M)$. 
\end{example}

\begin{sit}
    Let us recall the meaning and structure of the maps $f$
    and $g$ for such an Auslander context, cf.
    \cite[II.4]{Bas}. The $\End_{A}(M)$--bimodule
    homomorphism $f:M\otimes_{A}M^{*}\to \End_{A}(M)$ is
    given by $f(m,\lambda)(m')=m\cdot \lambda(m')$. It is
    the {\em norm map\/} of the module $M$, whose image
    consists of all endomorphisms that factor through a
    finite free, equivalently a finite projective,
    $A$-module. Its cokernel, isomorphic to the Morita
    defect $\oC(A,M)$, is by definition the {\em stable
    endomorphism ring\/} $\uEnd_{A}(M)$ of $M$ over $A$, and
    fits into the exact sequence of $\End_{A}(M)$-bimodules
    $$
    0\to \Omega_{\End_{A}(M\oplus A)/A} \to
    M\otimes_{A}M^{*}\xto{f}\End_{A}(M)\to \uEnd_{A}(M)\to
    0\,.
    $$
    The morphism of $A$-bimodules
    $g:M^{*}\otimes_{\End_{A}(M)}M \to A$ is simply the
    evaluation map, $g(\lambda, m) = \lambda(m)$, and its
    image is, by definition, the (twosided) {\em trace
    ideal\/} $\tau_{M}(A)\subseteq A$ of $M$ in $A$. Note
    that the trace ideal is all of $A$, equivalently, $g$ is
    bijective, if and only if $M$ is an
    $A$--generator.
\end{sit}

In other words, the vanishing of the Morita defects for an
Auslander context has the following classical
interpretation, see \cite[II.Prop.4.4]{Bas}:

\begin{proposition}
\label{prop:ME}
    Let $C=C(A,M)$ be an Auslander context as just described.
    \begin{enumerate}
    \item $\oC = \uEnd_{A}(M)$ vanishes if and only if $M$
    is a finite projective $A$-module.
    
    \item $\oC' = A/\tau_{M}(A)$ vanishes if and only if $M$
    is an $A$-generator.
    \qed
    \end{enumerate}
\end{proposition}

\begin{sit}
\label{sit:AC}
    Every Morita context $(C,e)$ defines several Auslander
    contexts. For example, $Ce$ is a right $A=eCe$-module
    that splits as $Ce = e'Ce\oplus eCe = M\oplus A$. The
    left $C$-module structure on the $(C,A)$-bimodule $Ce$
    provides a ringhomomorphism
    \begin{diagram}
	\alpha_{C}:C &\rTo^{\quad} &
	\Hom_{eCe}(Ce,Ce)=\End_{A}(M\oplus A)\,,
    \end{diagram}
    in detail,
    \begin{diagram}
	C=\begin{pmatrix}
	    B&M\\
	    N&A
	\end{pmatrix}&\rTo^{\alpha_{C}=
	\begin{pmatrix}
	    \beta&\id_{M}\\
	    g^{\sharp}&\id_{A}
	\end{pmatrix}
	}&
	\begin{pmatrix}
	    \End_{A}(M)& M\\
            M^{*} & A
	\end{pmatrix}=\End_{eCe}(Ce)\,,
    \end{diagram}
    where $\beta:B\to\End_{A}(M)$ defines the $B$-module
    structure on the $(B,A)$-bimodule $M$ and $g^{\sharp}:
    N\to M^{*}=\Hom_{A}(M,A)$ is adjoint to
    $g:N\otimes_{B}M\to A$. The idempotent $e\in C$,
    corresponding to $1_{A}$, maps to the idempotent
    corresponding to $\id_{A}$ in $\End_{A}(M\oplus A)$. 
\end{sit}

\begin{defn}
    We will say that $C$ {\em represents a {\em(\/}right{\em
    )\/} Auslander context\/}, on the pair $(A,M)$, if
    $\alpha_{C}$ is an isomorphism. In that case, we call
    $e\in C$ an {\em Auslander idempotent.\/}
\end{defn}

Auslander contexts can be characterized through the grade of
the Morita defect $\oC$, as we now show.

\begin{proposition}
\label{prop:AC}
    The canonical morphism of Morita contexts
    $$
    \alpha_{C}:C\to \Hom_{eCe}(Ce,Ce)=\End_{A}(M\oplus 
    A)
    $$
    \begin{enumerate}
    \item  is {\em injective\/} iff $\Hom_{C}(\oC,C)=0$ iff
    $\grade_{C}\oC\ge 1$;
    
    \item is {\em bijective\/} iff $\Ext^{i}_{C}(\oC,C)=0$ for
    $i=0,1$ iff $\grade_{C}\oC\ge 2$.
    \end{enumerate}
    In particular, $(C,e)$ is a (right) Auslander context if and 
    only if $\grade_{C}\oC \ge 2$.
\end{proposition}

\begin{proof}
    The proof consists essentially of applying
    $\Hom_{C}(-,C)$ to the exact sequence
    \ref{lem:fundamental}(\ref{eq:1}). Indeed, remark first
    that with $A=eCe$ one has
    \begin{align*}
	\Hom_{C}(Ce\otimes_{A}eC,C) &\cong
	\Hom_{A}(Ce,\Hom_{C}(eC,C)) \quad\text{by
	adjunction,}\\
        &\cong \Hom_{A}(Ce,Ce) \quad\quad\text{as 
        $\Hom_{C}(eC,C)\cong Ce$.}
    \end{align*}
    Moreover, with respect to these isomorphisms,
    $$
    \alpha_{C}\cong\Hom_{C}(\mu_{e},C):C\cong
    \Hom_{C}(C,C)\to \Hom_{A}(Ce,Ce) \cong
    \Hom_{C}(Ce\otimes_{A}eC,C)\,.
    $$
    Now split the exact sequence (\ref{eq:1}) into two short
    exact sequences,
    \diagramstyle[lefteqno]
    \begin{diagram}[eqno=(\dagger)]
    \label{eq:d}
    0&\rTo& CeC &\rTo &C &\rTo &\oC &\rTo &0\,,\\
    0&\rTo& \Omega_{C/A} &\rTo &Ce\otimes_{A}eC &\rTo &CeC
    &\rTo &0\,,
    \end{diagram}
    and apply $\Hom_{C}(-,C)$ to obtain the diagram with
    exact row and column
    \begin{diagram}[eqno=(\dagger\dagger)]
    \label{eq:dd}
    &&&&&&0\\
    &&&&&&\dTo\\
    0&\rTo &\Hom_{C}(\oC,C) & \rTo & C & \rTo &
    \Hom_{C}(CeC,C) & \rTo & \Ext^{1}_{C}(\oC,C) &\rTo 0\\
    &&&&&\rdTo_{\alpha_{C}\cong \Hom_{C}(\mu_{e},C)} &\dTo\\
    &&&&&& \Hom_{C}(Ce\otimes_{A}eC,C)\\
    &&&&&&\dTo\\
    &&&&&&\Hom_{C}(\Omega_{C/A},C)&.
    \end{diagram}
    The Ker-Coker-Lemma for the factorization of
    $\alpha_{C}$ displayed in this diagram shows immediately
    that $\alpha_{C}$ is injective iff $\Hom_{C}(\oC,C)=0$,
    and that is the first claim. Moreover, if
    $\Hom_{C}(\oC,C) = 0$, then also
    $\Hom_{C}(\Omega_{C/A},C) = 0$, by \ref{lem:grade}, as
    $(\Omega_{C/A})e = 0$ by \ref{lem:fundamental}. The
    second claim follows then again from the
    Ker-Coker-Lemma.
\end{proof}

\begin{corollary}
\label{MoritaAuslander}
For a Morita context $(C,e)$, the 
following are equivalent:
\begin{enumerate}
    \item
    \label{MoritaAuslander.1}
    $e$ is a Morita idempotent;

    \item
    \label{MoritaAuslander.2}  
    $M$ is $A$--projective, $\beta:B\to \End_{A}(M)$ 
    is an isomorphism, and $g^{\sharp}:N\to \Hom_{A}(M,A)$ 
    is an isomorphism.

    \item
    \label{MoritaAuslander.3}  
    $N$ is $A$--projective and $\lgrade_{C}\oC \ge 2$.
\end{enumerate}
\end{corollary}

\begin{proof} If $e$ is a Morita idempotent, then $(C,e)$ is
both a left and a right Auslander context by \ref{prop:AC},
and then \ref{prop:ME} shows both that
(\ref{MoritaAuslander.1}) $\Rightarrow$
(\ref{MoritaAuslander.2}) and that (\ref{MoritaAuslander.1})
$\Rightarrow$ (\ref{MoritaAuslander.3}). The isomorphisms in
(\ref{MoritaAuslander.2}) exhibit $C$ as a right Auslander
context on $(A,M)$, and the implication
(\ref{MoritaAuslander.2}) $\Rightarrow$
(\ref{MoritaAuslander.1}) is then \ref{prop:ME} again.
Finally, (\ref{MoritaAuslander.3}) simply states that aside
from $N$ being $A$--projective, $C$ is isomorphic to the
Auslander context $C(A^{\op},N)$ by the natural map $C\to
\End_{eCe^{\op}}(Ce)=\End_{A^{\op}}(N\oplus A)^{\op}$. It is
thus dual to (\ref{MoritaAuslander.2}), hence equivalent to
(\ref{MoritaAuslander.1}) as the latter statement is self
dual.
\end{proof}

Next, we will show that the grade of the Morita defect for
an Auslander context measures the vanishing of
$\Ext^{i}_{A}(M\oplus A,M\oplus A)$. To this end, we use
first the exact sequences ($\dagger$) in the preceding proof to
establish the following auxiliary results:

\begin{lemma}
\label{lem:aux}
If $\grade_{C}\oC = g \ge 2$, and $D$ is a finite projective
$C$--module, then
\begin{equation}
\tag{a}
\label{eq:2}
    \Ext^{i}_{C}(CeC,D) \cong
    \begin{cases}
        0 &\text{for $1 \le i < g-1$,}\\
    \Ext^{g}_{C}(\oC,D) &\text{for $i = g-1$.}      
    \end{cases}
\end{equation}
\begin{equation}
\tag{b}
\label{eq:3}
    \Ext^{i}_{C}(CeC,D)
    \xto{\cong} 
    \Ext^{i}_{C}(Ce\otimes_{A}eC,D)\quad\text{for $0\le i \le 
    g-1$.}
\end{equation}
In particular, $\Ext^{i}_{C}(Ce\otimes_{A}eC,D)\cong
\Ext^{i+1}_{C}(\oC,D)$ for $1\le i \le g-1$.
\end{lemma}

\begin{proof} For each $i\ge 1$, the first exact sequence in
($\dagger$) above yields isomorphisms
$$
\Ext^{i}_{C}(CeC,D)\xto{\cong} \Ext^{i+1}_{C}(\oC,D)\,,
$$
whereas the second one yields exact sequences
\begin{diagram}
\Ext^{i-1}_{C}(\Omega_{C/A},D)&\rTo &
\Ext^{i}_{C}(CeC,D) & \rTo & \Ext^{i}_{C}(Ce\otimes_{A}eC,D) &
\rTo &\Ext^{i}_{C}(\Omega_{C/A},D)
\end{diagram}
Finally, use that $\Ext^{i}_{C}(\Omega_{C/A},D) = 0$ for
$i<g$, by \ref{lem:grade}.
\end{proof}

We also need the following.

\begin{lemma}
\label{lem:tor}
For each $j\ge 1$, the $C$--bimodule $\Tor_{j}^{eCe}(Ce,eC)$
is annihilated by $e$ on either side, thus naturally a
$\oC$--bimodule. Furthermore, $\Tor_{j}^{eCe}(Ce,eC) \cong
\Tor_{j}^{eCe}(e'Ce,eCe')= \Tor_{j}^{A}(M,N)$.
\end{lemma}

{\sc Proof}.
Multiplying $eC$ with $e$ from the right factors as
\begin{xxalignat}{3}
    &\hphantom{\square}
	    &(-)e&\colon eC\onto eCe \into eC
		    &&\square
			    \end{xxalignat}

Now we can establish the significance of the grade of the
Morita defect for an Auslander context. Indeed, we have the 
following, more general result.

\begin{theorem}
\label{thm:grade}
    Let $A$ be a ring and $M$ an $A$--generator. The grade
    $g$ of $\uEnd_{A}(M)$ as $\End_{A}(M)$--module
    satisfies
    $$
    g = 1 + \inf\{1 \le i \le \infty \mid
    \Ext^{i}_{A}(M,M) \ne 0\}\ge 2\,.
    $$    
    Moreover, if $g$ is finite, then
    $$
    \Ext^{g}_{\End_{A}(M)}(\uEnd_{A}(M),\End_{A}(M)) \cong
    \Ext^{g-1}_{A}(M,M)\,.
    $$
\end{theorem}

\begin{proof} Let $C =\End_{A}(M\oplus A)$ be the Auslander
context defined by $A$ and $M$, and note that $e'C$ is a
finite projective $C$--module. Now consider the two spectral
sequences with the same limit,
\begin{align*}
'E^{i,j}_{2}
&=\Ext^{i}_{C}(\Tor_{j}^{A}(Ce,eC),e'C)\Longrightarrow 
\bbbe^{i+j}\\
''E^{i,j}_{2}
&=\Ext^{i}_{A}(Ce,\Ext^{j}_{C}(eC,e'C))\Longrightarrow 
\bbbe^{i+j}\,.
\end{align*}
The second one degenerates, with $''E^{i,j}_{2}=0$ for
$j\neq 0$, as $eC$ is a projective $C$--module. Moreover,
$$
\Hom_{C}(eC,e'C)\cong e'Ce\cong M
$$ 
as $A$--modules, whence the limit term satisfies
$$
\bbbe^{*} \cong \Ext^{*}_{A}(Ce,e'Ce)\cong  
\Ext^{*}_{A}(A\oplus M,M)\,.
$$
By \ref{lem:tor}, the terms $\Tor_{j}^{A}(Ce,eC)$ are
$\oC$--modules for $j\ge 1$, hence $'E^{i,j}_{2}$ vanishes
for $j\ge 1$ and $i< g$ by \ref{lem:grade}. Accordingly, the
edge homomorphisms $'E^{i,0}_{2}\to \bbbe^{i}$ are
isomorphisms for $i\le g$, thus,
$$
'E^{i,0}_{2} = \Ext^{i}_{C}(Ce\otimes_{A}eC,e'C) \cong 
\Ext^{i}_{A}(Ce,e'Ce)\quad \text{for $i\le g$.}
$$
As $C$ is an Auslander context, $g\ge 2$ by \ref{prop:AC},
and so 
$$ 
\Ext^{i}_{C}(Ce\otimes_{A}eC,e'C) \cong
\Ext^{i+1}_{C}(\oC, e'C) \quad \text{for $1\le i\le g-1$,}
$$ 
by \ref{lem:aux}. Putting these isomorphisms together, we
find for $1\le i<g-1$ that
$$
0=\Ext^{i+1}_{C}(\oC, e'C) \cong \Ext^{i}_{A}(Ce,e'Ce) =
\Ext^{i}_{A}(M\oplus A,M) \cong \Ext^{i}_{A}(M,M)\,,
$$
whereas for $i=g-1$, which is at least equal to $1$, we get
$$
\Ext^{g}_{C}(\oC, e'C)\cong \Ext^{g-1}_{A}(Ce,e'Ce) \cong
\Ext^{g-1}_{A}(M,M)\,.
$$
To complete the proof, it remains to show that
$$
\Ext^{g}_{C}(\oC, e'C) \cong
\Ext^{g}_{\End_{A}(M)}(\uEnd_{A}(M),\End_{A}(M))\,.
$$
As $M$ is a generator, $\oC'=0$ by \ref{prop:ME}, and so $B$
and $C$ are Morita equivalent in view of \ref{cor:ME}. An
explicit equivalence is given by $-\otimes_{C}Ce':\Mod C\to
\Mod B$, see \ref{prop:catME}, thus,
$$
\Ext^{g}_{C}(\oC, e'C) \cong \Ext^{g}_{B}(\oC\otimes_{C}Ce', 
e'C\otimes_{C}Ce')\,.
$$
To evaluate the right hand side, observe, as in the proof of
\ref{lem:gradechange}, that
$$
\oC\otimes_{C}Ce' = \oC e' =\oC \cong \uEnd_{A}(M)
$$ 
as $B$--(bi)module and that $e'C\otimes_{C}Ce' \cong e'Ce' =
B$. Thus, we get finally
$$
\Ext^{g}_{\End_{A}(M)}(\uEnd_{A}(M),\End_{A}(M))\cong 
\Ext^{g}_{C}(\oC, e'C) \cong \Ext^{g-1}_{A}(M,M)
$$
and the theorem is established.
\end{proof}

As a first application, we have the following.

\begin{corollary} 
\label{cor:strange}
Assume $M$ is an $A$--generator such that $\End_{A}(M)$ is a
right noetherian ring of finite global dimension, say,
$\gldim \End_{A}(M) = d$. If $\Ext^{i}_{A}(M,M) = 0$ for
$i=1,\ldots,d-1$, then $M$ is finite projective over $A$,
and $A$ is of global dimension equal to $d$.
\end{corollary}

\begin{proof} 
As $B=\End_{A}(M)$ is assumed to be right noetherian of
global dimension $d$, for any finitely generated $B$--module
$Y$, vanishing of $\Ext^{i}_{B}(Y,B)$ for $i\le d$ implies
$Y=0$. Applied to $Y = \uEnd_{A}(M)$, the preceding theorem
shows then that $\Ext^{i}_{A}(M,M) = 0$ for $i=1,\ldots,d-1$
implies $\uEnd_{A}(M) =0$. However, this means that both
Morita defects for $C = \End_{A}(M\oplus A)$ vanish, whence
that $M$ is a finite projective generator and that
$\End_{A}(M)$ is Morita equivalent to $A$.
\end{proof}

Theorem \ref{thm:grade} is formulated directly for the ring
$A$ and its generator $M$, with the associated Auslander
context hidden in the proof.  To display that context more
prominently, we make the following definition that is
motivated by \cite{Aus1, ARe}.

\begin{defn} A Morita context $(C,e)$ is a {\em Wedderburn 
context\/} if the following conditions are satisfied.
\begin{enumerate}
    \item  The Morita defect $\oC'$ vanishes.

    \item  The Morita defect $\oC$ satisfies 
    $\grade_{C}\oC \ge 2$.
\end{enumerate}
\end{defn}

Recall, \cite[Sect.8]{Aus1}, that a module $P$ over a ring
$R$ is a {\em Wedderburn projective\/} if it is finite
projective and the natural ring homomorphism
$$
R\lto \End_{\End_{R}(P)}(\Hom_{R}(P,R))
$$
is an isomorphism. With this notion, we can complement
\ref{thm:grade} as follows.

\begin{prop}
\label{prop:Wed}
For a Morita context $(C,e)$ the following 
conditions are equivalent.
\begin{enumerate}
    \item
    \label{prop:Wed.1}
    $(C,e)$ is a Wedderburn context.

    \item
    \label{prop:Wed.2}
    $(C,e)$ is a (right) Auslander context with 
    $M=e'Ce$ an $A=eCe$--generator.

    \item
    \label{prop:Wed.3}
    $(C,e')$ is a (right) Auslander context with $N=eCe'$ a
    Wedderburn projective over $B=e'Ce'$. Note that from
    this point of view, $\oC\cong B/\tau_{N}(B)$, the
    quotient of $B$ modulo the trace ideal of $N$.
\end{enumerate}
\end{prop}

\begin{proof} 
The equivalence (\ref{prop:Wed.1}) $\Longleftrightarrow$
(\ref{prop:Wed.2}) follows from \ref{prop:AC} and
\ref{prop:ME}. 

If $(C,e)$ is a Wedderburn context, then both $(C,e)$ and
$(C,e')$ are Auslander contexts by \ref{prop:AC}, and $N$ is
finite projective over $B$ by \ref{prop:ME}. That $(C,e')$
is an Auslander context implies $M\cong \Hom_{B}(N,B)$ as
$(B,A)$--modules and $A\cong \End_{B}(N)$ as rings, see
\ref{sit:AC}. That $(C,e)$ is an Auslander context yields,
by \ref{sit:AC} again, the first isomorphism in the chain
$$
B\cong \End_{A}(M) \cong 
\End_{\End_{B}(N)}(\Hom_{B}(N,B))\,,
$$
whereas the second simply substitutes what we just stated.
Thus, $N$ is a Wedderburn projective over $B$ and
(\ref{prop:Wed.1}) $\Longrightarrow$ (\ref{prop:Wed.3}) is
established. 

Finally, we show (\ref{prop:Wed.3}) $\Longrightarrow$
(\ref{prop:Wed.2}). As $N$ is $B$--projective and $(C,e')$
is isomorphic to the Auslander context $C(B,N)$ by
assumption, $\oC'=0$ by \ref{prop:ME}. Moreover, for any
projective $B$--module $N$, the $\End_{B}(N)$--module
$\Hom_{B}(N,B)$ is a generator, as is well known or follows
also directly from \ref{prop:ME} applied to the Morita
context $C(B,N)$. As $(C,e')$ is an Auslander context, the
natural ring homomorphism $A\to \End_{B}(N)$ is an
isomorphism, and the implication follows.

For the final comment in (\ref{prop:Wed.3}), observe that
$\oC\cong \Coker(f: M\otimes_{A}N\to B)$ by
\ref{lem:fundamental}, and that $f$ can be identified with
the evaluation map $\Hom_{B}(N,B)\otimes_{\End_{B}(N)} N \to
B$, with image $\tau_{N}(B)$, as $(C,e')$ is an Auslander
context.
\end{proof}

\begin{rem}
    The classical duality between {\sc (Rings, Generators)}
    and\linebreak {\sc (Rings, Wedderburn Projectives)}, see
    \cite{Aus1,ARe}, appears here simply as the exchange
    $e\rightleftarrows e'$ of the complementary idempotents
    in the Wedderburn context $(C,e)$.
\end{rem}

\section{Equivalences to the Generalized Nakayama
Conjecture}

The result in \ref{thm:grade} allows to shed further light
onto the relation among various homological conjectures that
have been originally formulated for Artin algebras and their
modules; see \cite{Yam} for a detailed account of how little
is yet known about these conjectures!

Auslander and Reiten always maintained that many of these
conjectures should indeed be true for rings that are finite
as modules over a noetherian commutative ring. In this
section only, we call such rings {\em Noetherian algebras}.

For example, they stated the following form of the
Generalized Nakayama Conjecture, first for Artin algebras in
\cite[p.70]{ARe}, later in the more general context here,
see \cite[Introduction to Chap.V]{ASW}:

\begin{Conjecture}[GNC] 
Let $A$ be a Noetherian algebra. If $M$ is a finitely
generated $A$--generator, then
$$
\Ext^{i}_{A}(M,M) = 0\quad\text{for $i>0$}\quad
\Longrightarrow \quad\text{$M$ is projective.}
$$
\end{Conjecture}

For Artin algebras, it was shown in the same paper
\cite{ARe} that this conjecture is equivalent to:

\begin{Conjecture}[GNC'] 
For any simple module $S$ over an Artinian algebra $A$,
there is an integer $n\ge 0$ such that
$\Ext^{n}_{A}(S,A)\neq 0$.
\end{Conjecture}

This conjecture was strengthened for Artin Algebras by Colby
and Fuller in \cite{CFu} to the Strong Nakayama Conjecture
that we formulate here again for Noetherian algebras:

\begin{Conjecture}[SNC] 
Let $A$ be a Noetherian algebra. If $M$ is a finitely
generated $A$--module, then
$$
\Ext^{i}_{A}(M,A) = 0\quad\text{for $i\ge
0$}\quad\Longrightarrow \quad\text{$M = 0$.}
$$
\end{Conjecture}

Obviously, (SNC) $\Rightarrow$ (GNC') for Artin algebras. In
general, (SNC) implies immediately what one might call the
Idempotent Nakayama Conjecture:

\begin{Conjecture}[INC] 
If $(C,e)$ is a Wedderburn context with $C$ a Noetherian
algebra, then
$$
\Ext^{i}_{C}(\oC,C) = 0\quad\text{for $i\ge
0$}\quad\Longrightarrow \quad\text{$\oC = 0$,}
$$
in other words, infinite grade of $\oC$ should force a
Wedderburn context to be a Morita equivalence.
\end{Conjecture}

Indeed, the following formulation that we will show to be
equivalent to {\sc (INC)} in a moment, indicates that this
conjecture is, on the face of it, considerably weaker than
{\sc (SNC)}, if only for the reason that over many rings a
Wedderburn projective module is automatically a generator.

\begin{Conjecture}[INC'] 
If $N$ is a Wedderburn projective over a Noetherian
algebra $B$, then
$$
\Ext^{i}_{B}(B/\tau_{N}(B),B) = 0\quad\text{for $i\ge
0$}\quad\Longrightarrow \quad\text{$N$ is a generator.}
$$
\end{Conjecture}

As a consequence of \ref{thm:grade}, we now show that also 
for Noetherian algebras the Strong Nakayama Conjecture 
implies the Generalized one, and that the two conjectures 
mentioned last are equivalent to {\sc (GNC)}.

\begin{prop} For Noetherian algebras and finitely generated 
modules over them, one has
$$
\text{\sc (SNC)} \Longrightarrow \text{\sc
(INC')}\Longleftrightarrow \text{\sc
(INC)}\Longleftrightarrow \text{\sc (GNC)}\,.
$$
\end{prop}

\begin{proof} As already remarked, the first implication is
obvious. For the equivalences, we only have to note that we
stay within the class {\sc (Noetherian Algebras, Finitely
Generated Modules)} when we pass from a Morita context to
its components, or associate to such a pair its Auslander
context. Indeed, $C$ is a Noetherian algebra if and only if
$A, B$ are Noetherian algebras and $M,N$ are finitely
generated modules over each $A$ and $B$, see
\cite[1.1.7]{MRo}, and if $A$ is a Noetherian algebra and
$M$ a finitely generated $A$--module, then the Auslander
context $\End_{A}(M\oplus A)$ is clearly as well a
Noetherian algebra. The equivalences $\text{\sc (INC')}
\Longleftrightarrow\text{\sc (INC)} \Longleftrightarrow
\text{\sc (GNC)}$ then follow from \ref{thm:grade} and
\ref{prop:Wed}.
\end{proof}

\begin{remarks}
(1) The conjectures {\sc (INC')}, {\sc (INC)}
    trivially hold if the algebras $C, B$ involved are
    already commutative noetherian rings. However, there
    seems to be no real advantage gained in either {\sc
    (SNC)} or {\sc (GNC)} if one assumes that $A$ is already
    commutative. In this sense, the aforementioned
    conjectures truly belong to the realm of (slightly)
    noncommutative algebra.

(2) One could as well state the above conjectures for
    arbitrary rings and modules. However, R.~Schulz
    \cite{Sch} gave a counterexample to {\sc (GNC)}, even
    for a self injective artinian ring $A$. Of course, his
    construction relies heavily upon the fact that the ring
    in question is not finite over its centre.

(3) The Generalized Nakayama Conjecture for Artin
    Algebras implies the so-called Tachikawa conjectures for
    Artin Algebras; see again \cite{Yam}. For a treatment of
    those conjectures for commutative noetherian rings, see
    \cite{ABS}.
\end{remarks}

\section{The Comparison Homomorphism for Hochschild Cohomology}

Let $(C,e)$ be a Morita context and $(A,B,M,N)$ the
associated Pierce components of $C$.  Additionally, we
assume form now on that $C$ comes equipped with a
$K$--algebra structure over some commutative ring $K$.  The
rings $A,B$ inherit then a $K$--algebra structure from $C$,
and $M,N$ become symmetric $K$--modules.

The aim here is to show that there exists a canonical
homomorphism of rings 
\[
\chi_{C/A}:\Hoch(C)\lto \Hoch(A)\,.
\]
The existence of such a homomorphism was already observed in
\cite{GSo}, though under several additional assumptions.

Before introducing the homomorphism, we first recall two
general facts from homological algebra that will be used 
below.

\begin{sit}
    All unadorned tensor products, enveloping algebras, Bar
    resolutions, Hoch\-schild cohomology groups, and the
    like are from now on taken over $K$.  For a $K$--algebra
    $A$, its enveloping algebra is $A^{\ev}=A^{\op}\otimes
    A$, so that an $A$--bimodule $\calm$ whose underlying
    $K$--module is symmetric is the same as a right
    $A^{\ev}$--module through $m(a^{\op}\otimes a') = ama'$. 
    Let $\bbbb(A)$ denote the Bar resolution or acyclic
    Hochschild complex of $A$ over $K$, and set
    \[
	\Hoch(A,\calm) = H(\Hom_{A^{\ev}}(\bbbb(A),\calm))\,,
	\]
    the Hochschild cohomology of $A$ over $K$ with values in
    an $A$--bimodule, more precisely an $A^{\ev}$--module,
    $\calm$.  We abbreviate further and set $\Hoch(A) =
    \Hoch(A,A)$.
\end{sit}

Note the following general fact on the multiplicative struture of 
Hochschild cohomology.

\begin{lemma}
    \label{lem:generalbar}
Let $R$ be a $K$--algebra.  If $X\to Y$ is a morphism of
complexes of $R^{\ev}$--modules that is a homotopy
equivalence if considered as a morphism of complexes of
$K$--modules, then the induced map of complexes
$$
\Hom_{R^{\ev}}(\bbbb R,X)\to \Hom_{R^{\ev}}(\bbbb R,Y)
$$
is a quasiisomorphism.
\end{lemma}

\begin{proof}
    Taking the mapping cone over $X\to Y$, it suffices to
    show that $\Hom_{R^{\ev}}(\bbbb R,Z)$ is acyclic as soon
    as $Z$ is a complex of $R^{\ev}$--modules that is
    contractible as complex of $K$--modules.  Viewing
    $\Hom_{R^{\ev}}(\bbbb R,Z)$ as the total (product)
    complex associated to the double complex
    $\Hom_{R^{\ev}}(\bbbb R_{n},Z_{m})$, each ``column'',
    that is, the complex for a fixed $n$, is acyclic, even
    contractible, as the terms $\bbbb R_{n}$ are induced from
    $K$--modules.  As $\bbbb R_{n} = 0$ for $n<0$, the
    double complex is concentrated in the ``right''
    halfplane, whence the total complex remains acyclic.
\end{proof}

In view of the fact just recalled, the augmentation of the Bar
resolution yields an isomorphism in cohomology
$$
H(\Hom_{R^{\ev}}(\bbbb R,\bbbb R))\lto \Hoch(R)\,,
$$
and the composition or Yoneda product on
$\Hom_{R^{\ev}}(\bbbb R,\bbbb R)$ induces the natural ring
structure on $\Hoch(R)$.  It is the same as the one
originally defined by M.~Gerstenhaber \cite{Ger}.
   
We will also use below the following basic property of the
Yoneda product.

\begin{lemma}
    \label{lem:quasialg}
Let $f:X\to Y$ be a morphism of complexes over some abelian
category and assume that the map $\Hom(X,f)$ is a
quasiisomorphism. The composition
$$
H(\Hom(X,f))^{-1}\circ H(\Hom(f,Y):H(\Hom(X,X))\to H(\Hom(Y,Y))
$$
is then a ring homomorphism with respect to Yoneda product on 
source and target.
\end{lemma}

\begin{proof}
    Indeed, the hypothesis simply means that for any morphism 
    of complexes $g:Y\to Y$, there exists a morphism of complexes
    $h:X\to X$ that renders the following square commutative
    \begin{diagram}
        X&\rTo^{f}& Y\\
	\dDashto^{h}&&\dTo_{g}\\
	X&\rTo^{f}& Y
    \end{diagram}
    and that such $h$ is unique up to homotopy.  The
    indicated map in cohomology sends the cohomology class
    of $g$ to that of $h$, and uniqueness up to homotopy
    then guarantees that the image of a composition is
    homologous to the composition of the images.
\end{proof}

Now we return to Morita contexts, with notation as 
introduced before.

\begin{prop} 
\label{lem:incl}
Let $(C,e)$ be a Morita context.  There is a canonical
semi-split inclusion of complexes of $C$--bimodules
\begin{diagram}
  {\tilde \mu_{e}}: Ce\otimes_{A}\bbbb(A)\otimes_{A}eC &\rInto&
  \bbbb(C)
\end{diagram}
that lifts the multiplication map $\mu_{e}:Ce\otimes_{A}eC
\to C$.

The adjoint to the inclusion ${\tilde \mu_{e}}$, given by
\begin{diagram}
   {\tilde \mu_{e}^\sharp}: \bbbb(A)&\rInto&
   \Hom_{C^{\ev}}(Ce\otimes eC, \bbbb(C))\cong
   e\bbbb(C)e\,,
\end{diagram}
is a semi-split inclusion and quasiisomorphism of complexes
of $A^{\ev}$--modules that induces $\id_{A}$ in homology. 
Moreover, ${\tilde \mu_{e}^\sharp}$ is a homotopy equivalence of 
complexes of right, as well as of left $A$--modules.
\end{prop}

\begin{proof}
Clearly, 
$$
Ce\otimes_{A}A^{\otimes (*+2)}\otimes_{A}eC \cong Ce\otimes
(eCe)^{\otimes *}\otimes eC
$$
is a direct summand of the $C$--bimodule $C^{\otimes *+2}$,
and the differentials are compatible. Moreover, the induced 
map in the zeroth homology is clearly $\mu_{e}$.

As for the second assertion, $Ce\otimes eC \cong
(e^{\op}\otimes e )C^{\ev}$ is a finite, even cyclic,
projective $C^{\ev}$--module, whence
$$
   \Hom_{C^{\ev}}(Ce\otimes eC, \bbbb(C))\cong
   \bbbb(C)\otimes_{C^{\ev}}C^{\ev}(e^{\op}\otimes e)\cong 
   e\bbbb(C)e
$$
remains exact, resolving $C \otimes_{C^{\ev}}
C^{\ev}(e^{\op}\otimes e) \cong eCe = A$.  Furthermore,
${\tilde \mu_{e}^\sharp}$ realizes $A^{\otimes (*+2)} =
(eCe)^{\otimes (*+2)}$ as the obvious direct summand of
$eC\otimes C^{\otimes *}\otimes Ce$.  Direct inspection
shows that with respect to these identifications,
${\tilde \mu_{e}^\sharp}$ induces $\id_{A}$ in the only nontrivial
homology group.

As for the final claim, we note that the right $A$--linear
map sending an element $w$ to $e\otimes w\in eC\otimes
C^{\otimes *}\otimes Ce$ contracts the augmented complex
$e\bbbb(C)e\to A$.  The other side is left to the reader. 
Both augmented complexes $\bbbb A\to A$ and $e\bbbb Ce\to A$
are thus contractible as complexes of one-sided modules,
whence ${\tilde \mu_{e}^\sharp}$, lying over the identity on $A$,
must be a homotopy equivalence.
\end{proof}

\begin{sit}
    \label{sit:relativeHH}
    For every $C$--bimodule\footnote{tacitly assumed, as always, to be
	a symmetric $K$--bimodule} $\calm$, we set 
	\begin{gather*}
	\bbbb(C/A) := \Coker({\tilde \mu_{e}}:
	Ce\otimes_{A}\bbbb(A)\otimes_{A}eC \into \bbbb(C))
	\\
	\Hoch(C/A,\calm) := H(\Hom_{C^{\ev}}(\bbbb(C/A),\calm))
	\end{gather*}
    to obtain first a
    semi-split exact sequence of $C^{\ev}$--modules
    \begin{equation}
    \label{eq:33}
    \tag{*}
       \begin{diagram}
            0&\rTo& Ce\otimes_{A}\bbbb(A)\otimes_{A}eC &\rTo& 
           \bbbb(C) & \rTo & \bbbb(C/A) & \rTo & 0
       \end{diagram}
    \end{equation}
    and then long exact cohomology sequences
    $$
    \cdots\to \Hoch^{i}(C/A,\calm)\to \Hoch^{i}(C,\calm)\to
    \Hoch^{i}(A,e\calm e)\to \Hoch^{i+1}(C/A,\calm)\to
    \cdots\,,
    $$
    where we have used that
    $$
    \Hom_{C^{\ev}}(Ce\otimes_{A}\bbbb(A)\otimes_{A}eC,\calm) 
    \cong \Hom_{A^{\ev}}(\bbbb(A), e\calm e) 
    $$
    by adjunction.
\end{sit}

\begin{sit}
    \label{sit:comparison}
    If we take $\calm = C$, then $e\calm e = A$ and we get a
    natural map
    $$
    \chi^{*}_{C/A}=\Hoch^{*}(\mu_{e},C):\Hoch^{*}(C)\lto
    \Hoch^{*}(A)
    $$
    relating the Hochschild cohomology of $A$ to that of
    $C$.  According to the classical theorem of
    M.~Gerstenhaber \cite{Ger}, Hochschild cohomology of any
    algebra carries naturally the structure of a graded
    commutative algebra --- and much more structure besides. 
    We now show that $\chi^{*}_{C/A}$ respects at least the
    product structure.
\end{sit}

\begin{theorem}
    \label{thm:ringhomo}
    The map $\chi^{*}_{C/A}$ is a homomorphism of graded
    $K$--algebras. 
\end{theorem}

\begin{proof}
In view of our preparations, the claim follows essentially
from the commutativity of the following diagram:
\begin{diagram}
\Hom_{C^{\ev}}(\bbbb C,\bbbb C)&\rTo^{-\circ
{\tilde \mu_{e}}}&\Hom_{C^{\ev}}(Ce\otimes_{A}\bbbb
A\otimes_{A}eC,\bbbb C)\\
\dTo^{e(-)e}&&\dTo^{\cong}_{\text{adjunction}}\\
\Hom_{A^{\ev}}(e\bbbb Ce,e\bbbb Ce)&\rOnto^{-\circ {\tilde
\mu_{e}^\sharp}}&\Hom_{A^{\ev}}(\bbbb A,e\bbbb Ce)\\
&&\uTo_{{\tilde \mu_{e}^\sharp}\circ -}^{\simeq}\\
&&\Hom_{A^{\ev}}(\bbbb A,\bbbb A)
\end{diagram}
where $\cong$ denotes an isomorphism of complexes, $\simeq$
a quasiisomorphism.  Namely, commutativity of the square is
straightforward from the definition of $\mu_{e}$ and the
fact that $Ce\otimes_{A}(-)\otimes_{A}eC$ is left adjoint to
$e(-)e\cong eC\otimes_{C}(-)\otimes_{C}Ce$.  The composed
map along the top and down the right side induces by
definition $\chi^{*}_{C/A}$ in cohomology, taking into
account \ref{lem:generalbar}, once for $R=C$ and once for
$R=A$.  The vertical map at the left is visibly a ring
homomorphism, and ${\tilde \mu_{e}^\sharp}$ satisfies the
conditions in \ref{lem:quasialg}, in view of \ref{lem:incl}.
            
Taken together, we get 
$$
\chi^{*}_{C/A} = \left(H({\tilde \mu_{e}^\sharp}\circ
-)^{-1}\circ H(-\circ{\tilde \mu_{e}^\sharp})\right)\circ
H(e(-)e)\,,
$$
as a composition of two homomorphisms of graded rings.
\end{proof}

\section{Hochschild Cohomology of Auslander Contexts}

The starting point for the results in this section is the classical
observation that Hoch\-schild theory is Morita invariant. 
More precisely, by specializing the argument in
\cite[1.2.7]{Lod}, one has the following.

\begin{prop}[Morita Invariance of Hochschild Theory]
\label{prop:MI}
    If the Morita defect $\oC$ vanishes, then the inclusion
    $\mu_{e}$ of complexes of $C^{\ev}$--modules in {\em
    \ref{lem:incl}\/} is a homotopy equivalence,
    equivalently, $\bbbb(C/A)$ is a contractible complex of
    $C^{\ev}$--modules. In particular, for every 
    $C^{\ev}$--module $\calm$, the natural map 
    $$
    \Hoch^{i}(C,\calm)\to \Hoch^{i}(A,e\calm e)
    $$
    is an isomorphism.\qed
\end{prop}

\begin{sit}
    What can we say in case the Morita defect does not
    vanish?  This clearly comes down to understanding the
    cohomology groups $\Hoch^{*}(C/A,\calm)$ introduced in
    \ref{sit:relativeHH}.  To be able to control those
    groups through the hypercohomology spectral sequence,
    {\em we assume henceforth that $C$ is projective as
    $K$--module}.  Under this assumption, each of the direct
    $K$--summands $A,B,M,N$ of $C$ is a projective
    $K$--module as well, and the terms of the complexes in
    the exact sequence (\ref{eq:33}) from
    \ref{sit:relativeHH} are projective $C^{\ev}$--modules.
\end{sit}

Moreover, we can easily determine the homology groups of
those complexes.

\begin{lem}
\label{lem:homology}
The homology of 
$Ce\otimes_{A}\bbbb(A)\otimes_{A}eC$ satisfies
$$
H_{*}(Ce\otimes_{A}\bbbb(A)\otimes_{A}eC) = 
\Tor^{A}_{*}(Ce,eC)\,.
$$
The homology groups of $\bbbb(C/A)$ are the $\oC$--bimodules 
$$
\oH_{j}:= H_{j}(\bbbb(C/A)) =  
\begin{cases}
\oC &\text{if $j=0$,}\\
\Omega_{C/A}&\text{if $j=1$,}\\
\Tor_{j-1}^{A}(Ce,eC)\cong \Tor_{j-1}^{A}(M,N)&\text{if $j>1$.}
\end{cases}
$$
\end{lem}

\begin{proof}
    The first statement is classical: As $A$ and $Ce$ are
    projective $K$--modules, $Ce\otimes_{A}\bbbb(A)$
    constitutes a projective resolution of $Ce$ as
    $A$--module. The second statement follows from the first
    and the long exact homology sequence associated to
    (\ref{eq:33}), as $\bbbb(C)$ is a resolution of $C$, and
    the exact sequence of low order homology is canonically
    identified with the fundamental sequence in
    \ref{lem:fundamental},
\begin{diagram}
    0&\rTo& \Omega_{C/A}&\rTo& Ce\otimes_{A} 
    eC&\rTo^{\mu_{e}} & C&\rTo& \oC&\rTo& 0\\
    && \dTo^{\cong} && \dTo^{\cong} && \dTo^{\cong} &&
    \dTo^{\cong}\\
    0&\rTo& \oH_{1}&\rTo&
    H_{0}(Ce\otimes_{A}\bbbb(A)\otimes_{A}eC)
    &\rTo&H_{0}(\bbbb(C))&\rTo&
    \oH_{0}&\rTo& 0
\end{diagram}
Finally, \ref{lem:fundamental} shows that $\Omega_{C/A}$ is
naturally a $\oC$--bimodule, and \ref{lem:tor} shows the
same for $\Tor_{j}^{A}(M,N)$ when $j\ge 1$.
\end{proof}

As the complex $\bbbb(C/A)$ is bounded in the direction of
the differential and consists of projective
$C^{\ev}$--modules, the cohomology $\Hoch(C/A)$ of
$\Hom_{C^{\ev}}(\bbbb(C/A),C)$ can be calculated through the
hypercohomology spectral sequence
$$
E^{i,j}_{2}=\Ext^{i}_{C^{\ev}}(\oH_{j},C)\quad \Longrightarrow 
\quad \bbbe^{i+j}= \Hoch^{i+j}(C/A)\,.
$$
The next result evaluates the low order terms of this 
spectral sequence relative to the grade of the Morita defect.

\begin{lemma}
\label{lem:hc}
If $(C,e)$ is a Morita context with $\grade_{C}\oC = g$, 
then
$$
\Hoch^{i}(C/A) \cong
\begin{cases}
0 &\text{for $i < g$,}\\   
\Hoch^{0}(C,\Ext_{C}^{g}(\oC,C)) &\text{for $i = g$.}
\end{cases}
$$
\end{lemma}

{\sc Proof}.
First note that $\grade_{C}\oH_{j} \ge g$ by
\ref{lem:homology} and \ref{lem:grade}, and that thus 
$$
g = \inf\{j + \grade_{C}\oH_{j}\mid j\ge 0\}
$$ 
as $\oH_{0}=\oC$. We now analyze the $E_{2}$--terms in the
hypercohomology spectral sequence above using yet another
spectral sequence, see \cite[XVI.4(5)]{CEi}:

If $X, Y$ are $(R,S)$--bimodules over $K$--projective
$K$--algebras $R,S$, then each $\Ext^{q}_{S}(X,Y)$ is an
$R$--bimodule and there is a spectral sequence
$$
{}'E^{p,q}_{2}=\Hoch^{p}(R,\Ext^{q}_{S}(X,Y))\quad
\Longrightarrow \quad
{}'\bbbe^{p+q}=\Ext^{p+q}_{R^{\op}\otimes S}(X,Y)\,.
$$
Taking $R=S=C, X=\oH_{j}, Y=C$, we get the implications
\begin{alignat*}{2}
    &\Ext^{q}_{C}(\oH_{j},C) = 0&\quad &\text{for $q <
    \grade_{C}\oH_{j}$}\\
    \Longrightarrow\quad &\Ext^{q}_{C^{\ev}}(\oH_{j},C) =
    0&\quad &\text{for $q < \grade_{C}\oH_{j}$}\\
    \Longrightarrow\quad &\Hoch^{i}(C/A) = 0&\quad
    &\text{for $i < g =\inf\{ j+\grade_{C}\oH_{j}\mid j\ge
    0\}$.}
\end{alignat*}
As thus $E^{i,j}_{2}=0$ for $i<g$, the edge homomorphism
$E^{g,0}_{2}\to \bbbe^{g}$ is an isomorphism, whence 
$$
\Ext^{g}_{C^{\ev}}(\oC,C) = E^{g,0}_{2} \xto{\cong} \bbbe^{g}
= \Hoch^{g}(C/A)\,.
$$
On the other hand, ${}'E^{p,q}=0$ for $q < g$ shows that the
edge homomorphism ${}'\bbbe^{g}\to {}'E^{0,g}_{2}$ is an
isomorphism, whence for $\oH_{0}=\oC$ we get
\begin{xxalignat}{3}
    &\hphantom{\square}
	        &\Ext^{g}_{C^{\ev}}(\oC,C)&= {}'\bbbe^{g} \xto{\cong}
			{}'E^{0g}_{2} = \Hoch^{0}(C,\Ext_{C}^{g}(\oC,C))\,.
			            &&\square
						                \end{xxalignat}
This analysis of the relative Hochschild cohomology groups
$\Hoch^{*}(C/A)$ translates immediately into the following 
result.

\begin{theorem} 
\label{thm:MHH}
If $(C,e)$ is a Morita context with $\grade_{C}\oC =g$,
then the maps
$$
\chi^{j}_{C/A}:\Hoch^{j}(C)\lto \Hoch^{j}(A)
$$
are isomorphisms for $j \le g-2$, and there is an exact
sequence
$$
0\to \Hoch^{g-1}(C)\xto{\chi^{g-1}} \Hoch^{g-1}(A) \to
\Hoch^{0}(C,\Ext^{g}_{C}(\oC,C)) \to
\Hoch^{g}(C)\xto{\chi^{g}} \Hoch^{g}(A)\,.
$$
\qed
\end{theorem}

We list a few consequences.

\begin{cor} 
If $(C,e)$ is an Auslander context, then $C$ and $A$ have
isomorphic centres via $\chi^{0}_{C/A}:\calz(C)\xto{\cong}
\calz(A)$, and every outer derivation on $C$ is induced from 
a derivation on $A$.\qed
\end{cor}

\begin{cor}
    If the Morita defect $\oC$ is of infinite grade, then
    $\chi^{\hphantom *}_{C/A}:\Hoch(C)\to \Hoch(A)$ is an
    isomorphism.\qed
\end{cor}

We now give an application to deformation theory. Recall that if the 
algebra $A$ is projective over $K$, then its infinitesimal first 
order deformations over $K$ are parametrized by 
$\Hoch^{2}(A)$. Every formal deformation of $A$ over $K$ 
is trivial, that is, the algebra is {\em rigid\/}, if, and only 
if, $\Hoch^{2}(A) =0$. Similarly, if $M$ is an $A$--module, 
then $\Ext^{1}_{A}(M,M)$ describes the infinitesimal first 
order deformations of $M$ as $A$--module. Again, every 
formal deformation of $M$ as $A$--module is trivial, and $M$ 
is then called {\em rigid\/} as $A$--module, if, and 
only if, this extension group vanishes. In light of these 
facts, we have the following application.

\begin{cor}
    Let $A$ be a $K$--algebra and $M$ a module over it, with
    $A$ and $M$ projective as $K$--modules.
    If  $M$ is rigid as $A$--module and $A$ is rigid as
    $K$--algebra, then $\End_{A}(M\oplus A)$ is rigid too as 
    a $K$--algebra.
\end{cor}

\begin{proof}
    By \ref{thm:grade}, vanishing of $\Ext^{1}_{A}(M,M)$
    means that $\grade_{C}\oC \ge 3$ for the Auslander
    context $C=\End_{A}(M\oplus A)$.  Consequently,
    $\Hoch^{2}(C)$ embeds into $\Hoch^{2}(A)$ by
    \ref{thm:MHH}, and so $C$ is rigid along with $A$.
\end{proof}

Coming back to the general situation, in case the Morita
defect $\oC'$ vanishes, we get more precise information in
form of a spectral sequence that relates directly the
Hochschild cohomology of $B$ to that of $A$.

\begin{theorem}
\label{thm:SS}
    Let $(C,e)$ be a Morita context with $g=\grade_{C}\oC$.
    If $e'$ is a Morita idempotent, then there is a spectral
    sequence of graded $K$--algebras
    \begin{align*}
        E^{i,j}_{2}&=\Hoch^{i}(B,\Ext^{j}_{A}(M,M))\quad 
        \Longrightarrow\quad \bbbe^{i+j}=\Hoch^{i+j}(A)\,,
    \end{align*}
    and $E^{i,j}_{2}= 0$ for 
    $1\le j < g-1$. Moreover, there are natural $K$--algebra 
    homomorphisms
    $$
    \Hoch(B)\xto{\Hoch(B,\beta)}
    \Hoch(B,\End_{A}(M)) = 
    E^{\bdot,0}_{2}\xto{d}\bbbe=\Hoch(A)\,,
    $$
    induced by the ring homomorphism $\beta:B\to
    \End_{A}(M)$, as in \ref{sit:AC}, and the edge
    homomorphism $d$.  
    
    If furthermore $g\ge 2$, that is, if $(C,e)$ is a
    Wedderburn context, then $\Hoch^{i}(B)\cong E^{i,0}_{2}$
    for each $i\ge 0$, and the edge homomorphisms provide
    isomorphisms
    $$
    \Hoch^{i}(B) \cong \Hoch^{i}(A)\quad \text{for $i \le
    g-2$.}
    $$
    The exact sequence of low order terms becomes then
    $$
    0\to \Hoch^{g-1}(B)\to \Hoch^{g-1}(A) \to
    \Hoch^{0}(B,\Ext^{g-1}_{A}(M,M)) \to
    \Hoch^{g}(B)\to \Hoch^{g}(A)\,.
    $$
\end{theorem}

\begin{proof}
As $e'$ is a Morita idempotent, $M$ is a $B$--projective
$A$--generator and $A\cong \End_{B^{\op}}(M)^{\op}$ by
\ref{prop:Morita}.

Now any $(B,A)$--bimodule $M$ gives rise to two spectral
sequences of graded algebras with the same limit, see
\cite[XVI.4(5)]{CEi} again:
\begin{align*}
    {}'E^{i,j}_{2}=\Hoch^{i}(B,\Ext^{j}_{A}(M,M))\quad 
    &\Longrightarrow\quad \Ext^{i+j}_{B^{\op}\otimes 
    A}(M,M)\,,\\
    {}''E^{i,j}_{2}=\Hoch^{i}(A,\Ext^{j}_{B^{\op}}(M,M))\quad 
    &\Longrightarrow\quad \Ext^{i+j}_{B^{\op}\otimes 
    A}(M,M)\,.
\end{align*}
In our case, $M$ is projective as left $B$--module so that
$$
\Ext^{j}_{B^{\op}}(M,M) = 0\quad \text{for $j > 0$}
$$
and the second spectral sequence degenerates. Moreover,
$A\xto{\cong} \Hom_{B^{\op}}(M,M)$ as $A$--bimodules,
identifying the limit with $\Hoch(A)$. Combining this fact
with \ref{thm:grade}, the result follows.
\end{proof}

\begin{rem} As $B$ and $C$ are Morita equivalent when $e'$
is a Morita idempotent, $\Hoch(C)\cong \Hoch(B)$ in the
preceding result, and modulo this identification, the
isomorphisms and the exact sequence in \ref{thm:MHH}
are precisely the isomorphic edge homomorphisms and exact
sequence of low order terms resulting from the partial
degeneration of the spectral sequence above.

The spectral sequence in \ref{thm:SS} has already been used
in \cite{BLi} to investigate the Hochschild cohomology of
Artin algebras of finite representation type as well as the
behaviour of Hochschild cohomology under pseudo-tilting.
\end{rem}

\section{Invariant Rings: The General Case}

The aim of this section and the next is to apply our results
so far to invariant rings with respect to finite groups. For
background material on skew group rings see \cite{MRo} and
\cite{ARS}. We first recall the notions of a {\em
separable\/} ring homomorphism and of the {\em Noether} or
{\em homological} different, see \cite{AGo, HSu}.

\begin{defn} 
\label{defn:sep}
A homomorphism $f:R\to S$ of rings is {\em separable\/} if
the epimorphism of $S$--bimodules $\mu_{f}=\mu_{S/R}:
S\otimes_{R}S\to S, s'\otimes s''\mapsto s's'',$ splits.
\end{defn}

\begin{sit}
    A crucial measure of (non-)separability is the {\em
    Noether different\/}, defined by E.~Noether for
    homomorphisms $f$ between commutative rings and
    generalized as {\em homological different\/} by
    Auslander-Goldman in \cite{AGo} for associative
    algebras. To define it, recall that for a ring $S$ and
    an $S$--bimodule $\calm$, the $S$--invariants in $\calm$
    are given by the subgroup
    $$
    \calm^{S}=\{m\in \calm\mid ms=sm\quad\text{for all $s\in 
    S$}\}\,,
    $$
    a (symmetric) module over the centre $\calz(S) := S^{S}$ of
    $S$. If we view $S$ as an algebra over some commutative
    ring, say over $\bbbz$, then $\calm^{S}\cong
    \Hom_{S^{\ev}}(S,\calm)\cong \Hoch^{0}(S,\calm)$.
\end{sit}

\begin{defn}
Let $f:R\to S$ be a ring homomorphism. The {\em Noether 
different\/} of $S$ over $R$, or rather of $f$, is
$$
\theta(S/R) := \Imm\left(\mu^{S}_{S/R}:(S\otimes_{R}S)^{S}\to 
S^{S}\right)\subseteq \calz(S)\,,
$$
an ideal in the centre of $S$.
\end{defn}

The following two classical facts on the Noether different,
taken from \cite{AGo}, will be relevant in the next section:

\begin{sit}
    The ring $S$ is separable over $R$ if, and only if,
    $\theta(S/R) = \calz(S)$.
\end{sit}

\begin{sit}
\label{sit:diff}
    If $R$ is commutative noetherian, and if $f$ turns $S$ into
    a module--finite $R$--algebra, then for any prime ideal
    $\fp\subset R$, the localization $f_{\fp}:R_{\fp}\to
    S\otimes_{R}R_{\fp}$ is separable if, and only if,
    $\fp\not\subseteq f^{-1}(\theta(S/R))$. In other words,
    $V(f^{-1}(\theta(S/R)))\subseteq \Spec R$ is the locus over
    which $f$ is not separable.
\end{sit}

Now we turn our attention to group actions. Suppose given a
finite group $G$ acting through automorphisms on a ring $S$.
Denote by $SG = S\#G$ the corresponding skew
group ring, and by $R=S^{G}$ the ring of $G$--invariants in
$S$.

\begin{sit}
    Identifying the unit element $\epsilon\in G$ with the
    multiplicative identity in $SG$, the canonical ring
    homomorphism $S\to SG$ can be thought of as
    multiplication with $\epsilon$ on either right or left.
    
    Set $f=\sum_{g\in G}g\in SG$ and note that $fg=f=gf\in
    SG$ for every $g\in G$. In particular, $f^{2}=|G|f$,
    whence $f$ is not quite an idempotent, and the
    $SG$--modules $SGf$ or $fSG$ are not quite projective.
    However, the compositions
    \begin{diagram}
    S&\rTo^{(-)\epsilon}& SG &\rTo^{(-)f}&SGf \\
    S&\rTo^{\epsilon(-)}&SG &\rTo^{f(-)}&fSG 
    \end{diagram}
    are always bijections, thus, $SGf = Sf, fSG=fS$, and in
    this way $S$ can be viewed both as a left or right
    $SG$--module. The subring $R\subseteq S\subseteq SG$
    commutes with $f$, and so $S\cong Sf$ becomes an
    $(SG,R)$--bimodule, or, through $S\cong fS$, an
    $(R,SG)$--bimodule. The multiplication map on
    $fS\otimes_{SG}Sf$ takes its values in $R\subseteq SG$,
    and taken together, these data define a Morita context
    $$
    (C,e) =
    \begin{pmatrix}
        SG & Sf\\
        fS & R
    \end{pmatrix}\,.
    $$
\end{sit}

\begin{sit}
\label{sit:facts}
    The following facts are well known and easy to
    establish:
    \begin{itemize}
	\item The Morita defect $\oC'$ is isomorphic to
	$R/\tr_{G}(S)$, where
        $$
	\tr_{G}:S\to R\quad,\quad \tr_{G}(s)=\sum_{g}g(s)\,,
        $$ 
	is the $R$-linear trace map from $S$ to $R$. As
	$|G|r = \tr_{G}(r)$ for every $r\in R$, one has
	$|G|\oC' =0$.
	
	\item The Morita defect $\oC$ is isomorphic to
	$\oSG := SG/(SG)f(SG)=SG/SfS$.
        
	\item The Morita context $(C,e')$ is an Auslander
	context, that is,
        $$
	C\cong \End_{e'Ce'}(Ce')\cong \End_{SG}(SG\oplus
	fS)\,,
        $$
        or, in detail,
        $$
        C \cong
        \begin{pmatrix}
	 SG& Sf \cong \Hom_{SG}(fS,SG)\\
         fS & R\cong \End_{SG}(fS)
        \end{pmatrix}\,.
        $$	
    \end{itemize}
    In view of \ref{prop:AC}, one has thus always
    $\grade_{C}\oC' \ge 2$. The Morita defect $\oC'$
    vanishes if and only if the trace map $\tr_{G}$ is
    surjective, and then $fS$ is a projective $SG$--module.
    If these conditions are satisfied, then $SG$ and $C$ are
    Morita equivalent, see \ref{cor:ME}. In view of $|G|\oC
    =0$, these conditions are satisfied as soon as $|G|$ is
    invertible in $S$.
\end{sit}

To compare now $R$ with $C$ or $SG$, the crucial question
becomes: What is the grade of the Morita context $\oC$?
We first consider when $\oC$ vanishes, equivalently, when 
$(C,e')$ is a Wedderburn context.

\begin{lemma}
\label{lem:oC' vanishes}
If the Morita defect $\oC$ vanishes, then the inclusion
$R\into S$ is a {\em separable\/} ring homomorphism. 
\end{lemma}

\begin{proof} By \ref{lem:fundamental}, $\oC=0$ if and only if the 
restricted multiplication map 
\begin{equation}
\label{eq:4}
\tag{*}
    \mu_{f}:Sf\otimes_{R}fS = e'Ce\otimes_{eCe}eCe'\to e'Ce' = SG
\end{equation}
is a bijection, necessarily of $SG$--bimodules. The
$S$--bilinear map $\pi: SG\to S$ that sends any $g\in
G\setminus\{\epsilon\}$ to $0$ and $\epsilon$ to $1$ is
clearly a split epimorphism of $S$--bimodules. The
composition of these two maps is easily identified with the
multiplication map $\mu_{S/R}:S\otimes_{R}S\to S$.
\end{proof}

\begin{sit}
Before we investigate the converse of the preceding result,
we mention the following case straight out of
\cite[7.8]{MRo}: If $G$ acts through outer automorphisms on
$S$ and if $S$ is a {\em simple\/}\footnote{that is, there
are no nontrivial twosided ideals,} ring, then $SG$ is
simple too, whence $(SG)f(SG)=SG$, that is, $\oC=0$. This
leads to the following well known fact.
\end{sit}

\begin{prop}
    If $G$ acts on a simple ring $S$ through outer
    automorphisms, then $(C,e')$ is a Wedderburn context and
    $R$ is Morita equivalent to $C$. If, furthermore, the 
    trace map is surjective, then $R, C$ and $SG$ are Morita
    equivalent. If in the latter situation $C$ is a 
    $K$--algebra that is projective as $K$--module, then
    \begin{xxalignat}{3}
	&\hphantom{\square}
    &\Hoch(R)&\cong \Hoch(C) \cong \Hoch(SG)\,.
	&&\square
    \end{xxalignat}
\end{prop}
\begin{remark}
    This result has been used recently in \cite{AFLS} to
    determine the Hochschild cohomology of invariant
    subrings of the Weyl algebra with respect to finite
    group actions.
\end{remark}

We now exhibit a case, where separability of $S$ over $R$
forces vanishing of the Morita defect $\oC'$. This case will
cover most geometrically relevant situations. It is partly
motivated by Yoshino's treatment of quotient surface
singularities in \cite[10.8]{Yos}.

\begin{defn} 
A finite group $G$ acts {\em infinitesimally through outer
automorphisms\/} on $S$, if the natural inclusion of
$S$--bimodules $S\epsilon\into SG$ induces a bijection on
the subgroups of $S$--invariants,
$(-)\epsilon:S^{S}\xto{\cong} SG^{S}$. Equivalently, for
every $g\neq \epsilon$ and for every $s'\neq 0\in S$, there
exists an $s\in S$ such that
$$
ss'-s'g(s) \neq 0\,.
$$
\end{defn}

Note the following simple fact.

\begin{lemma}
If the finite group $G$ acts infinitesimally through outer
automorphisms on $S$, then the centre of $SG$ satisfies
$$
\calz(SG) \cong \calz(S)^{G}\,.
$$
In particular, if $S$ is commutative, then the centre of 
$SG$ equals $R$.
\end{lemma}

\begin{proof}
Calculate $(SG)^{SG}$ as $\left((SG)^{S}\right)^{G}$.    
\end{proof}

The following example provides a stock of such group
actions and explains their ubiquity in geometry.

\begin{example} 
\label{ex:dom}
If a finite group $G$ acts faithfully on a commutative
domain $S$, then it acts infinitesimally through outer
automorphisms. Indeed, $ss'-s'g(s) = (\epsilon-g)(s)\cdot
s' = 0$ for some $g\neq \epsilon$, some $s'\neq 0$, and all
$s\in S$ iff $\epsilon-g = 0$ on $S$ iff the action is not
faithful.
\end{example}

\begin{prop} 
Assume the finite group $G$ acts infinitesimally through
outer automorphisms on some ring $S$. The inclusion $R\into
S$ is then a separable ring homomorphism if and only if
$\oC=0$ if and only if $(C,e')$ is a Wedderburn context.
\end{prop}

\begin{proof} By \ref{lem:oC' vanishes}, we only need to
consider the case when $S$ is separable over $R$, whence
when there exists an $S$--bilinear splitting $\zeta:S\to
S\otimes_{R}S\cong Sf\otimes_{R}fS$ of the multiplication
map on $S$. Let $\tilde\zeta:S\to SG$ be the composition of
$\zeta$ with the multiplication map $\mu_{f}$ on $SG$ as in
(\ref{eq:4}). To establish that (\ref{eq:4}) is surjective,
it clearly suffices to show $\tilde\zeta(1) =1$. Indeed, on
the one hand, $s\tilde\zeta(1) = \tilde\zeta(s) =
\tilde\zeta(1)s$ for each $s\in S$, as $\tilde\zeta$ is by
construction $S$--bilinear. The assumption on the group
action then implies that $\tilde\zeta(1)\in S\epsilon\subset
SG$. On the other hand, if $\zeta(1)=\sum_{i}s_{i}'\otimes
s_{i}''$, then $\tilde\zeta(1) = \sum_{i,g}s_{i}'
g(s_{i}'')g$. Taken together, we get thus $\sum_{i}s_{i}'
g(s_{i}'') = 0$ for $g\neq \epsilon$, and so $\tilde\zeta(1)
= \sum_{i}s_{i}'s_{i}'' = 1$.
\end{proof}

Now we apply \ref{thm:SS} of the preceding section to obtain
the following result.

\begin{theorem}
    If $G$ acts through $K$--linear automorphisms on the
    $K$--projective $K$--algebra $S$, and if the trace map
    $\tr_{G}$ is surjective, then there is a spectral
    sequence 
    $$
    E^{i,j}_{2}=\Hoch^{i}(SG, \Ext^{j}_{R}(S,S))
    \Longrightarrow \Hoch^{i+j}(R)
    $$
    and natural homomorphisms $\Hoch^{i}(SG)\to
    \Hoch^{i}(SG,\End_{R}(S))= E^{i,0}_{2}\to \Hoch^{i}(R)$
    that are isomorphisms for $i\le \grade_{SG}\oSG -2$.
    \qed 
\end{theorem}

Moreover, it is known how to calculate the $E_{2}$--terms of
the preceding spectral sequence in terms of Hochschild
cohomology of $S$, see \cite{Lor, Ste}.

\begin{theorem}
\label{thm:deg}
    For any $SG$--bimodule $\calm$, there is a spectral sequence
    $$
    `E^{i,j}_{2}=H^{i}(G,\Hoch^{j}(S, \calm)) \Longrightarrow
    \Hoch^{i+j}(SG, \calm)\,.
    $$
    In particular, if $|G|$ is invertible in $K$, then
    $'E^{i,j}_{2}=0$ for $j\neq 0$, and the spectral sequence
    degenerates to yield
\begin{xxalignat}{3}
    &\hphantom{\square}
	    &\Hoch^{i}(SG, \calm)&\cong \Hoch^{i}(S, \calm)^{G}\,.
		    &&\square
			    \end{xxalignat}
\end{theorem}

In case $G$ acts infinitesimally through outer automorphisms, 
we can combine the two results into the following.

\begin{cor}
\label{cor:SSG}
If $G$ acts infinitesimally through outer $K$--linear
automorphisms on the $K$--projective $K$--algebra $S$, and
if $|G|$ is invertible in $K$, then there is a spectral 
sequence
$$
\Hoch^{i}(S, \Ext^{j}_{R}(S,S))^{G} \Longrightarrow
\Hoch^{i+j}(R)
$$
with natural homomorphisms
$$
\Hoch^{i}(S,SG)^{G}\cong \Hoch^{i}(S)^{G}\oplus \big(\bigoplus_{g\neq 
\epsilon}\Hoch^{i}(S, Sg)\big)^{G}\to\Hoch^{i}(R)\,,
$$
and these homomorphisms are isomorphisms for $i \le 
\grade_{SG}\oSG - 2$.\qed
\end{cor}

The remaining mysterious quantity is thus $\grade_{SG}\oSG$.
In case of commutative rings, it can be bounded from below
by the depth of the Noether different, the key point of the
next, and last, section.

\section{Invariant Rings: The Commutative Case}

>From now on, $S$ will be a commutative ring, and we will
assume that $G$ acts infinitesimally through outer
automorphisms, a condition that holds for faithful actions
on commutative domains by \ref{ex:dom}. In this case, the
Noether different of $S$ over its subring of invariants is
contained in the $S$--annihilator of the Morita defect
$\oC'\cong\oSG$ as we now show.

\begin{prop}
\label{prop:diff annihilates}
The kernel of the composed ring homomorphism $S\to SG \to
\oSG$ contains the Noether different $\theta(S/R)$, thus, 
equivalently, $\theta(S/R)\subseteq \Ann_{S}\oSG$.
\end{prop}

\begin{proof}
Set ${\mathfrak a} =\Ker(S\to \oSG)$ and note that ${\mathfrak a} =
S\cap \Ker(SG\to \oSG)$. With $X= S\times_{SG}
(Sf\otimes_{R}fS)$ the indicated fibre product, the image of
the first projection $p_{1}:X\to S$ is the ideal ${\mathfrak
a}$, as the rows in the commutative diagram
\begin{diagram}[textflow]
S\otimes_{R}S &\rTo^{\mu_{S/R}} & S\\
\uEqual&&\uTo^{\pi}\\
Sf\otimes_{R}fS &\rTo^{\mu_{SG}} & SG & \rTo & \oSG &\rTo 
&0\\
\uInto^{p_{2}}&&\uInto&&\uEqual\\
X &\rTo^{p_{1}} & S &\rTo & \oSG
\end{diagram}
of $S$--bimodules are exact. Note that $\pi(g) =
\delta_{g,\epsilon}$ and that $\pi$ retracts the inclusion
$S\subseteq SG$. As $S$ is commutative, $S=S^{S}$, and as
$G$ acts infinitesimally through outer automorphisms, $S =
S^{S}=(SG)^{S}$. Taking $S$--invariants preserves fibre
products, whence applying $(\ \;)^{S}$ to the first two
columns returns
\begin{diagram}
(S\otimes_{R}S)^{S} &\rTo^{\mu^{S}_{S/R}} & S\\
\uEqual&&\uEqual\\
(Sf\otimes_{R}fS)^{S} &\rTo^{\mu^{S}_{SG}} & (SG)^{S}\\
\uTo^{\cong}&&\uEqual\\
X^{S} &\rTo^{p^{S}_{1}} & S
\end{diagram}
By definition, $\mu^{S}_{S/R}$ has $\theta(S/R)$ as its
image, and the image of $p^{S}_{1}:X^{S}\to S$ is contained
in ${\mathfrak a}^{S} = {\mathfrak a}\subseteq S$.
\end{proof}

\begin{cor}
\label{cor:ineq}
If the ring $S$ is noetherian, then
$$
\grade_{SG}\oSG \ge \depth(\theta(S/R), S)\,.
$$
\end{cor}

\begin{proof} Deriving $\Hom_{SG}(X,Y) =
(\Hom_{S}(X,Y))^{G}$, for $SG$--bimodules $X,Y$, yields a
spectral sequence
$$
E^{i,j}_{2} = H^{i}(G,\Ext^{j}_{S}(X,Y)) \Longrightarrow
\Ext^{i+j}_{SG}(X,Y)\,.
$$
For $X=\oSG, Y=SG$, one obtains $\Ext_{S}(\oSG,SG) \cong
\oplus_{g}\Ext_{S}(\oSG,S)g$.  Now, as is well known,
$\Ext^{j}_{S}(\oSG,S) = 0$ for $j< \depth(\Ann_{S}\oSG, S)$,
see, e.g.~\cite[18.4]{Eis}, thus, $E^{i,j}_{2}=0$ for those
$j$.  Hence, $\Ext^{j}_{SG}(\oSG, SG)=0$ for
$j<\depth(\Ann_{S}\oSG, S)$, that is, $\grade_{SG}\oSG \ge
\depth(\Ann_{S}\oSG, S)$.  As we just saw,
$\theta(S/R)\subseteq \Ann_{S}\oSG$, whence
$\depth(\Ann_{S}\oSG, S)\ge \depth(\theta(S/R), S)$ and the
claim follows.
\end{proof}

We show next that the depth of $\theta(S/R)$ on $S$ controls
as well vanishing of the groups $\Hoch^{i}(S,Sg), g\neq \epsilon,$
that occurred in \ref{cor:SSG}.

\begin{lemma}
\label{lem:ann}
    Suppose that $S$ carries a $K$--algebra structure such
    that $G$ acts through $K$--linear isomorphisms.  In that
    case, $Sg$ is naturally an $S^{\ev}$--module for each
    $g\in G$, and there is an inclusion of $S$--ideals
    $$
    \theta(S/R)\subseteq \frac{\Ann_{S^{\ev}}(Sg) +
    \Ann_{S^{\ev}}(S)}{\Ann_{S^{\ev}}(S)} \subseteq
    \frac{S^{ev}}{\Ann_{S^{\ev}}(S)}= S\,.
    $$
    In particular, if $S^{\ev}$ is noetherian, then
    $\Ext^{i}_{S^{\ev}}(S,Sg) = 0$ for each $g\neq \epsilon$
    and $i < \depth(\theta(S/R), S)$.
\end{lemma}

\begin{proof} 
The proof is similar to that of \ref{prop:diff annihilates}.
First note that $K\subseteq R$ as $G$ acts through
$K$--linear automorphisms. Each $Sg$ is a cyclic
$S^{\ev}$--module, generated by $g$, and the resulting
epimorphism of $S$--bimodules $S^{\ev}\onto Sg$ factors
through the canonical surjective ring homomorphism
$S^{\ev}\onto S\otimes_{R}S$. In view of this, we may
replace $\Ann_{S^{ev}}$ in the claim everywhere by
$\Ann_{S\otimes_{R}S}$. Setting $I_{g} :=
\Ann_{S\otimes_{R}S}(Sg)$ and $i_{g}:I_{g}\into
S\otimes_{R}S$, we get the following diagram
\begin{diagram}[textflow]
&&&&0\\
&&&&\dTo\\
&&&&I_{\epsilon}\\
&&&&\dTo_{i_{\epsilon}}\\
0&\rTo&I_{g}&\rTo^{i_{g}}&S\otimes_{R}S&\rTo&Sg&\rTo&0\\
&&&\rdTo_{\mu_{S/R} i_{g}}&\dTo_{\mu_{S/R}}\\
&&&&S\\
&&&&\dTo\\
&&&&0
\end{diagram}
with exact row and column that shows the image of $\mu_{S/R}
i_{g}$ to be $(\Ann_{S^{\ev}}(Sg) +
\Ann_{S^{\ev}}(S))/{\Ann_{S^{\ev}}(S)}$. As $G$ acts
infinitesimally through outer automorphisms, $(Sg)^{S}=0$
for $g\neq \epsilon$, and so $ (I_{g})^{S}\cong
(S\otimes_{R}S)^{S} $ maps to
$\theta(S/R)=\Imm(\mu_{S/R}^{S})$ in $S$, whence the latter
ideal is contained in the image of $I_{g}$ in $S$.
\end{proof}

\begin{sit}
\label{sit:ass}
    Before formulating the key result in this section, we 
    summarize the necessary assumptions:
    \begin{itemize}
        \item  $S$ is a noetherian $K$--algebra that is 
        projective as a $K$--module and such that 
        $S^{\ev}$ is still a noetherian ring; for example, $S$ 
        is a finitely generated algebra over a field $K$;
    
        \item  the finite group $G$ acts infinitesimally through 
        outer $K$--linear automorphisms on $S$; for example,
        $S = K[V]/I$ is a domain, generated over $K$ by the 
        finite dimensional vector space $V$, and $G\subseteq 
        GL(V)$;
    
        \item  the invariant ring $R$ is again noetherian and 
        $S$ is a finitely generated $R$--module.
    
        \item  the order of $G$ is invertible in $K$.
    \end{itemize}
\end{sit}

Now we can state the main application of our
results to invariant rings. It is clearly the analogue of 
Schlessinger's theorem for tangent cohomology of invariant 
rings mentioned in the introduction.

\begin{theorem} 
\label{thm:key}
Under the assumptions of \ref{sit:ass}, one has
$$
g=\grade_{SG}{\oSG} \ge c=\depth(\theta(S/R),S)\,,
$$
and
$$
\Ext^{j}_{R}(S,S) = 0 \quad \text{for $1\le j \le g-2$.}
$$
There are natural isomorphisms
$$
\Hoch^{i}(S)^{G}\cong \Hoch^{i}(R) \quad\text{for $i \le
c-2$,}
$$
as well as an exact sequence
\begin{align*}
    0\to \Hoch^{c-1}(S)^{G}\to \Hoch^{c-1}(R)\to
    &\Hoch^{0}(SG,\Ext^{c-1}_{R}(S,S))\to \\
    &\to \Hoch^{c}(S,SG)^{G} \to \Hoch^{c}(R)\,.
\end{align*}
The groups $\Hoch^{i}(R)$ appear in general as limit of
a spectral sequence
$$
E^{i,j}_{2}=\Hoch^{i}(S,\Ext^{j}_{R}(S,S))^{G}\Longrightarrow
\Hoch^{i+j}(R)\,.
$$
\end{theorem}

\begin{proof} The inequality $c\le g$ is \ref{cor:ineq}. 
For the claim on $\Ext^{j}_{R}(S,S)$, we may assume that $g
\ge 2$, and in that case the Morita context $(C,e')$ from
\ref{sit:facts} is a Wedderburn context, as $|G|$ invertible
in $K$ forces $\oC'=0$.  The vanishing result then follows
from \ref{lem:gradechange} and \ref{thm:grade}.

Further, $\Hoch^{i}(S,Sg)\cong \Ext^{i}_{S^{\ev}}(S,Sg)$
for each $i$, as $S$ is projective over $K$, and those
groups vanish for $i < c$ by \ref{lem:ann}. Combining these 
facts, \ref{cor:SSG} now yields the remaining assertions.
\end{proof}

\begin{cor} Aside from the general assumptions in
\ref{sit:ass}, suppose furthermore that $G$ acts separably
outside isolated fixed points on a domain $S$ that is smooth
and finitely generated over a field $K$. In that case,
either $R$ is regular or
$$
c=\depth(\theta(S/R),S) = \dim S = g =\grade_{SG}\oSG\,,
$$
and one has
\begin{align*}
    \Hoch^{i}(R) &\cong \Hoch^{i}(S)^{G} \quad\text{for $i \le
    \dim S - 2$,}\\
    \Hoch^{\dim S -1}(R) &\cong \Hoch^{\dim S -1}(S)^{G}\oplus
    (\text{finite length $R$--module)}
\end{align*}
as $R$--modules. If $\dim S = 2$, the finite length module 
is zero.
\end{cor}

\begin{proof} 
As $S$ is supposed to be smooth over $K$, the global 
dimension of $SG$ equals the Krull dimension of $S$, and so
\ref{cor:strange} shows $g\le \dim S$, unless $R$ is regular.
The ring $S$ is in particular Cohen-Macaulay, 
and then $c = \depth(\theta(S/R),S) =\dim S$ by \ref{sit:diff},
as all minimal primes lying over $\theta(S/R)$ are by 
assumption maximal. Combined with the inequality $c\le g$, 
the claimed equalities follow. 

For the remainder, we only need to consider the case $g=c\ge
2$. As $S$ is smooth over $K$, its Hochschild cohomology
satisfies
$$
\Hoch^{i}(S) \cong \Hom_{S}(\Omega^{i}_{S/K},S)\,,
$$
by the Hochschild-Kostant-Rosenberg Theorem, see 
\cite[3.4.4]{Lod}, and so $\Hoch^{i}(S)$ is a projective 
$S$--module. In particular, its direct $R$--summand $\Hoch^{i}(S)^{G}$
is of maximal depth as $R$--module, whereas 
$$
\Ext^{\dim S -1}_{R}(S,S) \cong \Ext^{\dim S}_{SG}(\oSG,SG)
$$
is annihilated by $\theta(S/R)$ as $S$--module, and thus of
finite length. Accordingly, the image $X$ of $\Hoch^{\dim S
-1}(R)$ in $\Hoch^{0}(SG, \Ext^{\dim S -1}_{R}(S,S))$, is of 
finite length, and the resulting short exact sequence
$$
0\to \Hoch^{\dim S-1}(S,S)^{G}\to \Hoch^{\dim S -1}(R)\to X\to 0
$$
splits.

 Finally, note that for $\dim S = 2$ the
$R$--module $\Hoch^{1}(R) =\Der_{K}(R)$ is reflexive, thus
cannot contain a nontrivial submodule of finite length.
\end{proof}

We finish with the following example, due to M.~Auslander
\cite{Aus2}, that motivated our investigation of Morita contexts
and Hochschild cohomology.

\begin{example} Assume $S$ is the ring of formal 
power series in $n$ variables over a field $K$, and suppose
$G$ is a subgroup of $GL(n,K)$, with $|G|\in K^{*}$,
that acts through linear automorphisms on $S$. The group $G$ contains
no pseudo-reflections if, and only if, the branch locus
is of codimension at least 2. In that case, the natural map $SG\to
\End_{R}(S)$ is an isomorphism of rings, that is, $SG$ is an
Auslander context, as $S$ contains $R$ as an $R$--module
direct summand. If $n=2$, then $S$ is indeed a generator of 
the category of reflexive $R$--modules, and thus $SG$ is 
isomorphic to an Auslander algebra of the reflexive $R$--modules.
See, for example, \cite{Buc} for further information on this 
example and an investigation of how the Hochschild cohomology of 
the stable Auslander algebra or Morita defect $\oSG$ relates 
to that of $SG$.
\end{example}

\end{document}